\documentclass[review]{elsarticle}
\usepackage{tikz}
\usepackage{lineno,hyperref}
\usepackage{mathtools}

\pdfoutput=1
\modulolinenumbers[5]
\usepackage{float}
\RequirePackage{amsmath}[2016/11/05]
\usepackage[ruled,vlined]{algorithm2e}
\usepackage{natbib}
\journal{Engineering Applications of Artificial Intelligence}
\date{}
\usepackage{geometry}
\newgeometry{vmargin={15mm}, hmargin={12mm,17mm}}
\usepackage{caption}
\usepackage{subcaption}

\begin{document}

\begin{frontmatter}

\title{Improved Physics-informed neural networks loss function regularization with a variance-based term}

\author[1]{John M. Hanna}
\author[1,2]{Hugues Talbot}
\author[1]{Irene E. Vignon-Clementel}

\address[1]{Inria, Research Center Saclay Ile-de-France, France}
\address[2]{CentraleSupelec, Universit\'{e} Paris-Saclay}

\begin{abstract}

In machine learning and statistical modeling, the mean square or absolute error is commonly adopted as an error metric, also called a "loss function". While effective in reducing the average error, this approach may fail to address localized outliers, leading to significant inaccuracies in regions with sharp gradients or discontinuities. This issue is particularly evident in physics-informed neural networks (PINNs), where such localized errors are expected and affect the overall solution. To overcome this limitation, we propose a novel loss function that combines the mean and the standard deviation of the chosen error metric. By minimizing this combined loss function, the method ensures a more uniform error distribution and reduces the impact of localized high-error regions. The proposed loss function is easy to implement and tested on problems of varying complexity: 1D Poisson equation, unsteady Burger’s equation, 2D linear elastic solid mechanics, and 2D steady Navier Stokes equations. Results demonstrate improved solution quality and lower maximum error compared to the standard mean-based loss, with minimal impact on computational time.
\end{abstract}

\begin{keyword}
  Machine learning; outlier rejection; physics-informed neural networks; regularization; loss function
\end{keyword}

\end{frontmatter}


\section{Introduction}

Deep learning and its applications have been growing rapidly in the past decade or so. This is mainly due to the availability of large datasets and the groundbreaking technology of graphics processing units (GPUs) which have made it possible to train large neural network models with huge datasets. 

Recent neural network architectures have improved the accuracy of some deep learning tasks. For example, graph neural networks \cite{gnn} have made it easier to deal with data that can be organized in graphs, which has many applications such as in molecular biology \cite{gnn_molecular} or social networks \cite{gnn_social}. Transformer-based architectures \cite{trans_attension} with the novel self attention mechanism have significantly contributed to progress in natural language processing \cite{trans_nlp} and computer vision \cite{trans_vision} tasks.

In the recent years, neural networks were presented to solve partial differential equations (PDE), introducing physics-informed neural networks (PINN) \cite{pinn_main}. PINN reduce the need for big datasets and regularize outputs since the physics, governed by PDEs, are naturally regularizing the predictions. The main idea is to approximate the solution field of PDEs thanks to the strong approximation capabilities of neural networks \cite{universal_approx} . The solution is obtained by minimizing a combined loss function which represents the data mismatch (initial/ boundary conditions enforcement and/or sparse data in the domain) and the residual of the PDE over points sampled in the spatiotemporal domain called collocation points. The residual is usually obtained by automatic differentiation \cite{automatic_diff}. Since, the introduction of PINN, they have been applied in several fields including solid mechanics \cite{solid_pinn, pinn_hyperelastic}, fluid mechanics \cite{high_speed_pinn, pinn_fluid_review, pinn_hidden_fluid, pinn_arterial_flow}, flow in porous media \cite{adaptivity_hanna, pinn_porous, sensitivity_pinn}, topology optimization \cite{pinn_topology, pinn_topology_2}, composites manufacturing \cite{pinn_composite, self_supervised}, fracture mechanics \cite{pinn_fracture, pinn_fracture_2}, and many others.

One of the most crucial elements in deep learning, whether for supervised learning tasks or physics-informed neural networks tasks, is the choice of loss function. A pertinent loss function can lead to significant improvements in the optimization process, thus good accuracy can be reached more rapidly. In general, the loss function represents the mean error over the training points. 

In classical machine learning, many loss functions have been proposed, often associated with a specific task. The most common loss used in regression, the Mean Square Error (MSE) assumes that the error distribution is Gaussian with a uniform variance, which is often not realistic. In practice, the MSE works well enough even when the error distribution is only symmetric. However, when the variance of the error is not uniform (the heteroscedastic case), it is useful to use a loss that takes the local variance into account:

\[
\text{Loss}_{\text{H}} = \frac{1}{n}\sum_{i=1}^{n}\left( \frac{(y_i-\hat{y_i})^2}{2\hat{\sigma_i}^2}+\frac{1}{2}\log \hat{\sigma_i}^2\right),
\]

where $y_i$ is the true target value, $\hat{y_i}$ the predicted mean, $\hat{\sigma_i}^2$ the predicted variance. This model will predict both mean and variance around each training point, which is helpful when considering uncertainty. Such an approach was first proposed in risk-averse portfolio optimization in finance~\cite{markowitz1952portfolio}. In practice this loss is difficult to optimize~\cite{meanvariance_target_1994}. Gaussian processes~\cite{seeger2004gaussian}, in particular krieging~\cite{matheron1963principles} are also able to estimate localized mean and variance throughout the whole domain, but are not immediately applicable to neural networks, although there has been progress in this area~\cite{lee2017deep}. 

In these works, the main motivation is to estimate a local variance and adapt the MSE locally. However, dealing with non-uniform variance can also be seen as a form of robust regression, i.e. regression which is less sensitive to outliers. Loss functions like the mean absolute error or its differentiable variant, the Huber loss~\cite{huber1992robust} for regression tasks have been proposed in the literature to better deal with outliers. In practice, these losses are better suited for finding sparse solutions to regression problem~\cite{hastie2009elements}.

With PINN architectures, the situation is different. Classical outliers that represent unexpected out-of-distribution events can occur, but much more common are regions with high gradients or discontinuities that can lead to high localized residuals as compared to other regions, which in turn possibly deteriorates the predictions. These must be dealt with in a different fashion than in classical learning.

Several recent works attempt to improve the accuracy of PINN solution by adding extra regularization terms to the loss function. An example is the gradient-enhanced PINNs introduced by \cite{gpinn} where a higher order derivative term is added to the loss function. This term represents the residual gradient so as to ensure not just the residual of the PDE to be small, but also its derivative, which improves the solution accuracy. However it comes at a high computational cost, typically more than $(d+1)\times$ the original PINN computational time, where $d$ is the dimension of the problem.

Another work improved the training of PINN using a measure of the variance, or uncertainty of the output \cite{vpinn}. The main idea is to have two different outputs from the neural network, one representing the mean of the solution and the other represents the variance. In this sense it is similar to the heteroscedastic learning methods in classical learning mentioned above. Another term is added in the loss function that represents the log of the variance, the second output of the neural network. The idea is fundamentally different from our proposal in the sense that the variance of the solution is predicted by the neural network rather than being calculated. The idea can be rather seen as an automated point-wise scaling of the collocation points.

In this work, we present a new loss function that targets this issue. The core idea is based on an addition to the loss function. This addition includes the standard deviation of the chosen error to be minimized. By this addition, we are ensuring that minimizing the loss function can reduce both the mean and standard deviation of the error, thus encouraging smaller residuals than when minimizing only the mean. The method is tested on four PINN problems of different degrees of complexity: 1D Poisson equation, unsteady Burger's equation, 2D linear elasticity, and 2D steady-state Navier Stokes example.

The article is organized as follows: section~\ref{section:methods} is devoted to the introduction of the new proposed loss function, section~\ref{section:examples} offers different examples from physics-informed neural networks, and their results and comparison to standard loss function. Section~\ref{sec:discuss} finally  discusses several aspects such as computational costs or comparison to other well-known strategies, before drawing a conclusion on the proposed loss.

\section{Methodology}
\label{section:methods}

\subsection{Physics-informed neural networks}

Physics-informed neural networks was introduced by Raissi et al. \cite{pinn_main} as a method to solve forward and inverse problems governed by PDEs.

Given a PDE in the form of 

\begin{equation}
    u_t+\mathcal{N}(u;\lambda)=0, \ \ \ \ \ \mathbf{X}\in\Omega,\ t\in[0, T],
\end{equation}

where $u(t, \mathbf{X})$ is the solution to the PDE, $\mathcal{N}$ is a nonlinear differential operator parameterized by a material parameter $\lambda$, and $\Omega$ is the physical domain.

the residual of the PDE can be written as:

\begin{equation}
    r:=u_t+\mathcal{N}(u;\lambda).
\end{equation}

Initial and boundary conditions for the PDE are defined as

\begin{gather}
    u(0, \mathbf{X}) = u_0,\\
    u(t, \mathbf{X_D}) = u_D,\\
    \mathcal{B}(u(t, \mathbf{X_N})) = f(\mathbf{x_N}),
\end{gather}

where $\mathcal{B}$ is a differential operator, $\mathbf{X_D}$ is the boundary where Dirichlet boundary condition is defined, and $\mathbf{X_N}$ the boundary where the Neumann boundary condition is defined. 

The first step in solving this PDE with PINNs is to approximate the solution using a fully connected neural network. NNs are chosen due to their strong approximation capabilities \cite{universal_approx}.

The residual of the PDE is, then, calculated with automatic differentiation \cite{automatic_diff}. A combined loss function is formed with both the knowledge of the initial and boundary conditions as well as the physics represented by the residual of the PDE. The loss function is written as:

\begin{equation}
\mathcal{L} = \lambda_0\ \mathcal{L}_0 + \lambda_D\ \mathcal{L}_D + \lambda_N\ \mathcal{L}_N + \lambda_{f}\ \mathcal{L}_{r},
\end{equation}

where 

\begin{gather}
    \mathcal{L}_0 =\frac{1}{N_0}\sum_{i=1}^{N_0} ||u(t_0^i, \mathbf{X}_0^i) - u_0^i||^2,\\[2ex]
    \mathcal{L}_D =\frac{1}{N_D}\sum_{i=1}^{N_D} ||u(t_D^i, \mathbf{X}_D^i) - u_D^i||^2,\\[2ex]
    \mathcal{L}_N = \frac{1}{N_N}\sum_{i=1}^{N_N} ||B(u(t_N^i, \mathbf{X}_N^i))-f_N^i||^2,\\[2ex]
    \mathcal{L}_{r} = \frac{1}{N_r}\sum_{i=1}^{N_r} ||u_t+\mathcal{N}(u;\lambda)||_{(t_r^i, \mathbf{X}_r^i)}^2,
\end{gather}

where $\{t_0^i, \mathbf{x_0^i}, u_0^i\}_{i=1}^{N_0}$ are points for the initial condition, $\{t_D^i, \mathbf{x_D^i}, u_D^i\}_{i=1}^{N_D}$ the Dirichlet boundary condition, $\{t_N^i, \mathbf{x_N^i}, f_N^i\}_{i=1}^{N_N}$ the Neumann boundary condition, and $\{t_r^i, \mathbf{X_r^i}\}_{i=1}^{N_r}$ the collocation points where the physics is enforced.

A minimization problem is obtained which is commonly solved with gradient-based optimizers available in deep learning frameworks such as Adam \cite{adam} or a higher order methods \cite{bfgs}.

\subsection{Proposed loss function}

In general, common loss functions take the form of a mean of a chosen error:

\begin{equation}
    \mathcal{L} = \dfrac{1}{N}\sum_{i=1}^N e_i,
\end{equation}

where $N$ is the number of training points, $e_i$ is the chosen error function.

The $e_i$ term can represent the square error ${(\hat{y}_i(x) -y_i)}^2$, the absolute error ${|\hat{y}_i(x) -y_i)|}$, etc. where $\hat{y}$ is the approximation function and $y$ the label data.

The added term to the modified loss function represents the standard deviation of the error function $e$. The full loss function with the added term is written as:

\begin{equation}
    \mathcal{L} = \alpha \dfrac{1}{N}\sum_{i=1}^N e_i + (1-\alpha)\sqrt{\dfrac{\sum_{i=1}^N (e_i - \Bar{e})^2}{N}},
\end{equation}

where $\Bar{e}$ is the error mean which can be obtained as $\dfrac{\sum_{i=1}^N e_i}{N}$. $\alpha$ is a hyperparameter that takes values between $0$ and $1$ which identifies the contribution of each term in the loss function.

\section{Numerical examples}
\label{section:examples}

In this section, we show four problems solved with PINN. In all examples, a neural network architecture (NN) of 5 hidden layers and 20 neurons each is employed. The \texttt{tanh} activation function is employed for the hidden layers and linear function for the output layer. The Adam optimizer is chosen with a fixed learning rate of 0.001.

In all cases, we vary $\alpha$ from 0 to 1 to study the effect of this hyperparameter on the accuracy of the solutions, and we have the same initialization of neural network parameters (with Xavier initialization) each time $\alpha$ is changed to ensure a fair comparison. In all cases, NVIDIA RTX 3500 GPU was used for training.

\subsection{1D steady Poisson problem}

In this example, we solve the Poisson equation in one dimension as a simple test case. The strong form of the problem is written as:

\begin{equation}
    \begin{cases}
        -u_{xx} = 4x^2 \sin(x^2) - 2\cos(x^2), & x \in [-2\sqrt{\pi}, 2\sqrt{\pi}], \\
        u(-2\sqrt{\pi}) = u(2\sqrt{\pi}) = 1.
    \end{cases}
\end{equation}

The residual of the equation is defined as

\begin{equation}
    r(x) = u_{xx} + 4x^2 \sin(x^2) - 2\cos(x^2),
\end{equation}

and the solution field $u(x)$ is approximated with the neural network NN described above. The parameters of the neural network are optimized by minimizing a combined loss function that is written as:

\begin{equation}
    \mathcal{L} = \mathcal{L}_{\text{data}} + \mathcal{L}_{\text{PDE}},
\end{equation}
where for classical MSE training

\begin{equation}
    \mathcal{L}_{\text{data}} = \frac{1}{N_u} \sum_{i=1}^{N_u} \left| u_{\text{pred}}(x_u^i) - u^i_{\text{data}} \right|^2,
\end{equation}

and 

\begin{equation}
    \mathcal{L}_{\text{PDE}} = \frac{1}{N_f} \sum_{i=1}^{N_f} \left| r(x_f^i) \right|^2,
\end{equation}

while for the our new loss function

\begin{equation}
    \mathcal{L}_{\text{data}} = \frac{\alpha}{N_u} \sum_{i=1}^{N_u} \left| u_{\text{pred}}(x_u^i) - u^i_{\text{data}} \right|^2 + (1-\alpha)\sqrt{\dfrac{\sum_{i=1}^{N_u} (\left| u_{\text{pred}}(x_u^i) - u^i_{\text{data}} \right|^2 - \Bar{e})^2}{N_u}},
\end{equation}

where 

\begin{equation}
    \Bar{e} = \frac{1}{N_u} \sum_{i=1}^{N_u} \left| u_{\text{pred}}(x_u^i) - u^i_{\text{data}} \right|^2,
\end{equation}

and

\begin{equation}
    \mathcal{L}_{\text{PDE}} = \frac{\alpha}{N_f} \sum_{i=1}^{N_f} \left| r(x_f^i) \right|^2 + (1-\alpha)\sqrt{\dfrac{\sum_{i=1}^{N_f} (\left| r(x_f^i) \right|^2 - \Bar{e_r})^2}{N_f}},
\end{equation}

where

\begin{equation}
    \Bar{e_r} = \frac{1}{N_f} \sum_{i=1}^{N_f} \left| r(x_f^i) \right|^2,
\end{equation}

$N_u$ is the number of training points for the boundary conditions, $N_f$ is the number of collocation points for the PDE residual, $u_{\text{pred}}$ is the predicted solution from the neural network, $u_{\text{data}}$ is the data for boundary conditions, $x_u^i$ are the spatial coordinates for data points, $x_f^i$ are the spatial coordinates for collocation points.

\begin{figure}[H]
\centering
    \begin{subfigure}[b]{0.45\textwidth}
    \centering
    \includegraphics[width=250.0pt]{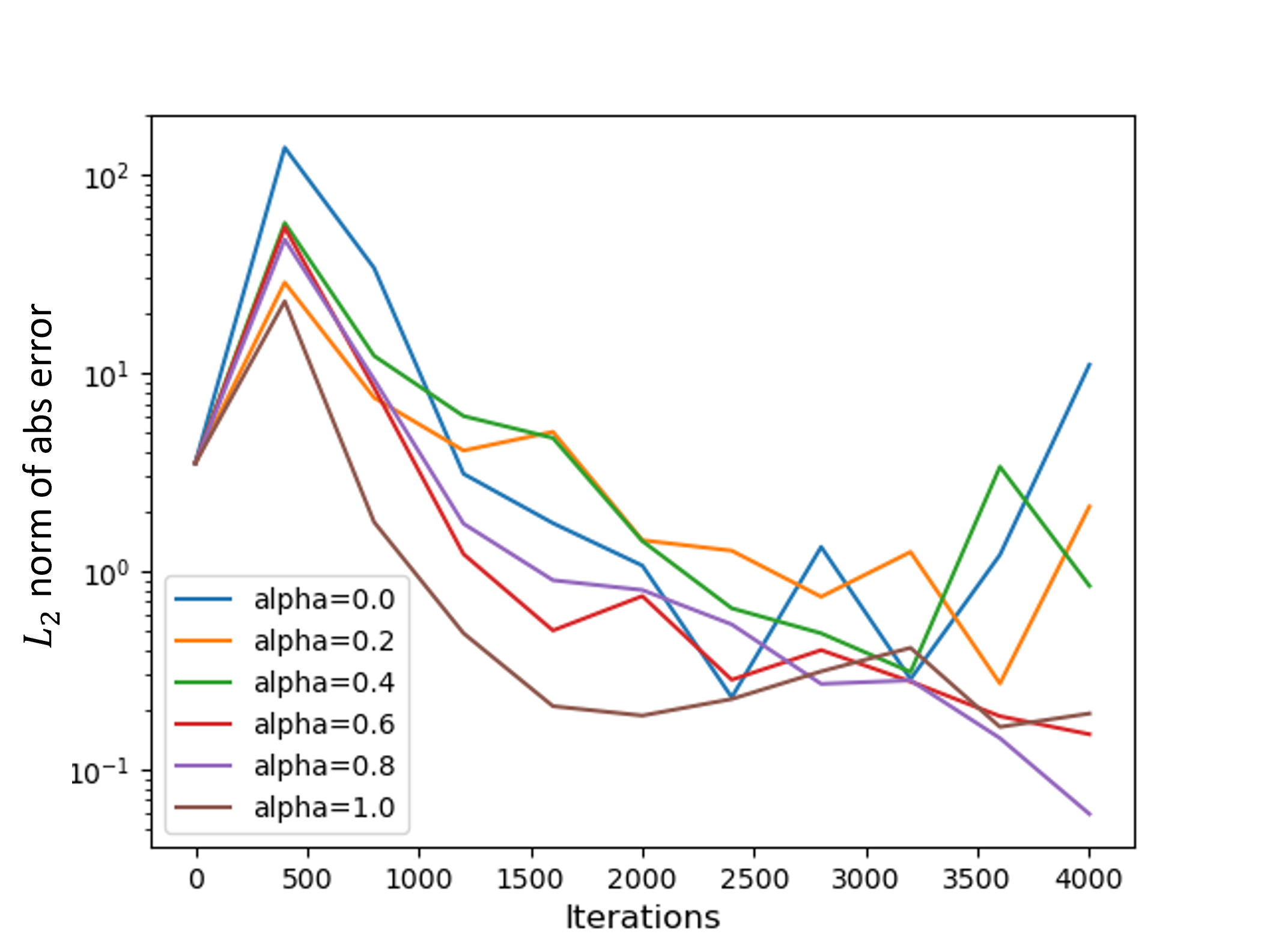}
    \end{subfigure}
    \begin{subfigure}[b]{0.45\textwidth}
    \centering
    \includegraphics[width=250.0pt]{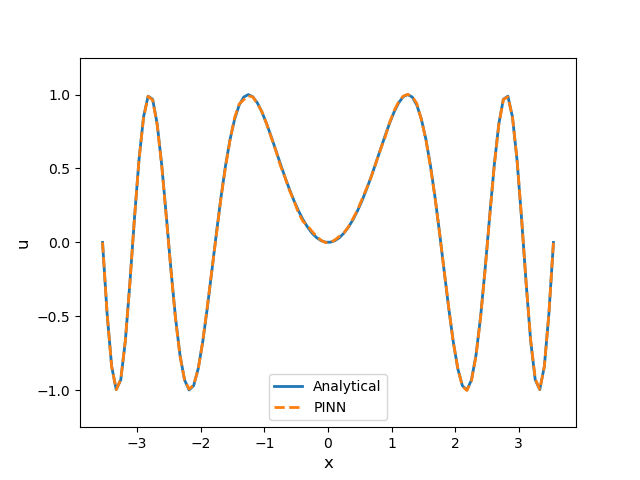}
    \end{subfigure}
\quad
\caption{On the left, the $L_2$ norm of the absolute error vs. the number of iterations for different values of $\alpha$ for the 1D Poisson example. On the right, a plot of the analytical solution along with the PINN solution for the best $\alpha$ (0.8).}
\label{fig:poissons_sols}
\end{figure}

The solutions are obtained with 4000 Adam iterations. Training is performed with 100 equally spaced collocation points, while the boundary conditions are exactly satisfied by construction. The computational cost using MSE loss is 8 seconds and 9 seconds when using the variance-based loss. Since an analytical solution to this problem exists, we plot the evolution of the $L_2$ norm of the absolute error between the PINN solution and the analytical solution, vs. the number of iterations for different values of $\alpha$, and the two solutions for the best $\alpha$ value in figure~\ref{fig:poissons_sols}.

As shown in figure~\ref{fig:poissons_sols}, the best solution is obtained with an $\alpha$ value of 0.8: adding the variance of the error does improve the accuracy of the prediction. However, the gain in accuracy is not large compared to the classical MSE loss ($\alpha=1$). This can be explained due to the simplicity and the smoothness of the solution. In the next examples, this improvement will be more pronounced.

\subsection{unsteady 1D Burgers' problem}

We now consider an unsteady example, which represents a standard test case in scientific machine learning, namely the unsteady Burgers' problem in one dimension as introduced by Raissi et al. \cite{pinn_main}. The Burgers' equation along with the initial and boundary conditions are written as:

\begin{equation}
    \begin{aligned}
        u_t + u u_x &= (0.01/\pi) u_{xx}, &\quad &x \in [-1, 1], \, t \in [0, 1], \\
        u(x, 0) &= -\sin(\pi x), &\quad &x \in [-1, 1], \\
        u(-1, t) &= u(1, t) = 0, &\quad &t \in [0, 1].
    \end{aligned}
\end{equation}

The residual of the Burger's equations is defined as

\begin{equation}
    r(x, t) = u_t + u u_x - (0.01/\pi) u_{xx},
\end{equation}

and the solution field $u(x, t)$ is approximated with the neural network NN. The parameters of the neural network are optimized by minimizing a combined loss function that is written as:

\begin{equation}
    \mathcal{L} = \mathcal{L}_{\text{data}} + \mathcal{L}_{\text{PDE}},
\end{equation}
where for classical MSE training

\begin{equation}
    \mathcal{L}_{\text{data}} = \frac{1}{N_u} \sum_{i=1}^{N_u} \left| u_{\text{pred}}(x_u^i, t_u^i) - u^i_{\text{data}} \right|^2,
\end{equation}

and 

\begin{equation}
    \mathcal{L}_{\text{PDE}} = \frac{1}{N_f} \sum_{i=1}^{N_f} \left| r(x_f^i, t_f^i) \right|^2,
\end{equation}

while for the our new loss function

\begin{equation}
    \mathcal{L}_{\text{data}} = \frac{\alpha}{N_u} \sum_{i=1}^{N_u} \left| u_{\text{pred}}(x_u^i, t_u^i) - u^i_{\text{data}} \right|^2 + (1-\alpha)\sqrt{\dfrac{\sum_{i=1}^{N_u} (\left| u_{\text{pred}}(x_u^i, t_u^i) - u^i_{\text{data}} \right|^2 - \Bar{e})^2}{N_u}},
\end{equation}

where 

\begin{equation}
    \Bar{e} = \frac{1}{N_u} \sum_{i=1}^{N_u} \left| u_{\text{pred}}(x_u^i, t_u^i) - u^i_{\text{data}}\right|^2,
\end{equation}

and 

\begin{equation}
    \mathcal{L}_{\text{PDE}} = \frac{\alpha}{N_f} \sum_{i=1}^{N_f} \left| r(x_f^i, t_f^i) \right|^2 + (1-\alpha)\sqrt{\dfrac{\sum_{i=1}^{N_f} (\left| r(x_f^i, t_f^i) \right|^2 - \Bar{e_r})^2}{N_f}},
\end{equation}

where

\begin{equation}
    \Bar{e_r} = \frac{1}{N_f} \sum_{i=1}^{N_f} \left| r(x_f^i, t_f^i) \right|^2,
\end{equation}

$N_u$ is the number of training points for the initial and boundary conditions, $N_f$ is the number of collocation points for the PDE residual, $u_{\text{pred}}$ is the predicted solution from the neural network, $u_{\text{data}}$ is the data for initial and boundary conditions, $x_u^i$ and $t_u^i$ are the spatial and temporal coordinates for data points, $x_f^i$ and $t_f^i$ are the spatial and temporal coordinates for collocation points.

\begin{figure}[H]
    \centering
    \includegraphics[width=0.7\linewidth]{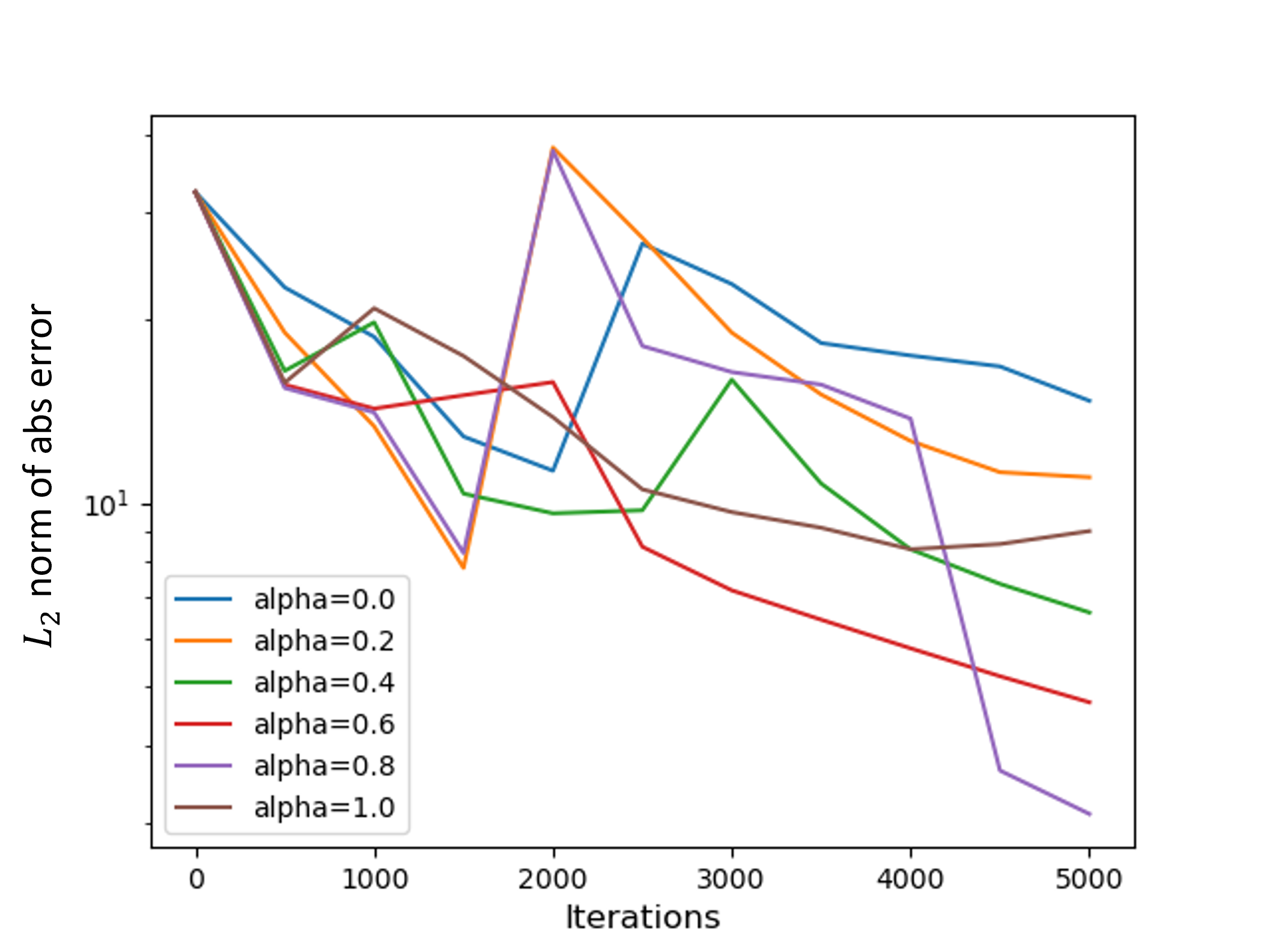}
    \caption{$L_2$ norm of the absolute solution error vs iterations for different values of $\alpha$ for the Burgers' problem.}
    \label{fig:burger_alphas}
\end{figure}

The solutions are obtained with 5000 Adam iterations. Training is carried out with 10,000 randomly sampled collocation points, while the initial and boundary conditions are automatically satisfied by construction. The computational cost using MSE loss is 10 seconds and 11 seconds when using the variance-based loss. For this case too, an analytical solution exists \cite{burgers_analytical}. First, the $L_2$ norm of the absolute error between the PINN and analytical solutions vs. iterations is plotted for different values of $\alpha$ as shown in figure~\ref{fig:burger_alphas}.

The effect of the added extra term is this time quite pronounced. It can be seen that solutions with $\alpha = 0.4, 0.6, 0.8$ perform better than the classical loss value ($\alpha = 1$), where the best solution is obtained with $\alpha=0.8$. 

\begin{figure}[H]
\centering
\begin{subfigure}{1\textwidth}
\centering
    \includegraphics[width=0.45\textwidth]{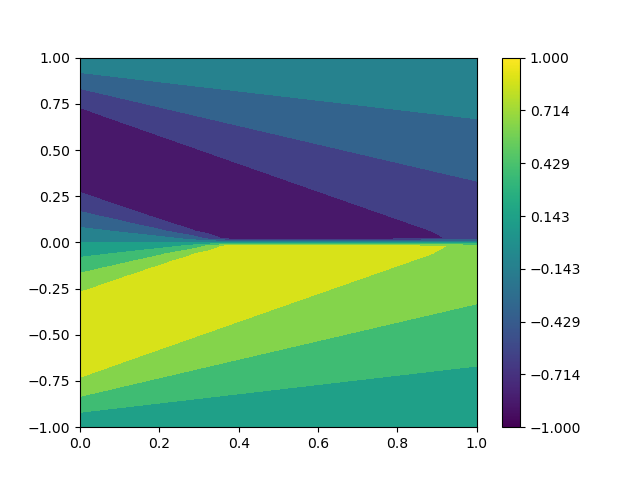}
    \caption{Analytical solution of the Burgers' problem.}
\end{subfigure}
\begin{subfigure}{0.45\textwidth}
    \includegraphics[width=\textwidth]{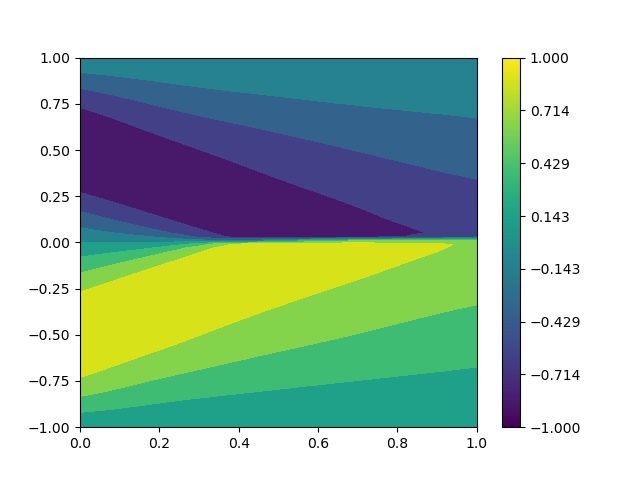}
    \caption{Burgers' solution using MSE loss with $\alpha=1$.}
\end{subfigure}
\begin{subfigure}{0.45\textwidth}
    \includegraphics[width=\textwidth]{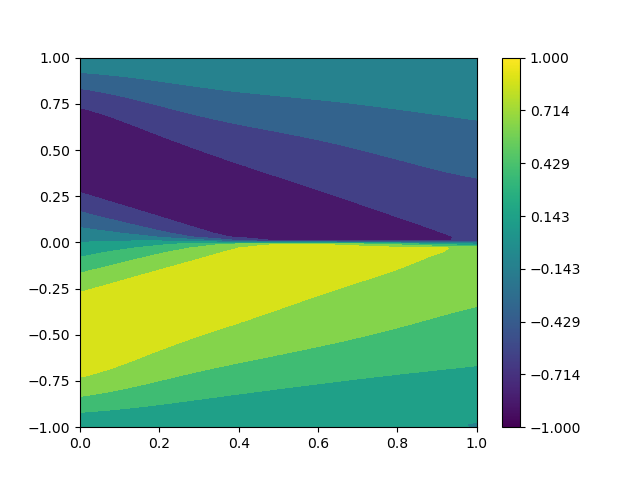}
    \caption{Burgers' solution using newly proposed loss with $\alpha=0.8$.}
\end{subfigure}
\begin{subfigure}{0.45\textwidth}
    \includegraphics[width=\textwidth]{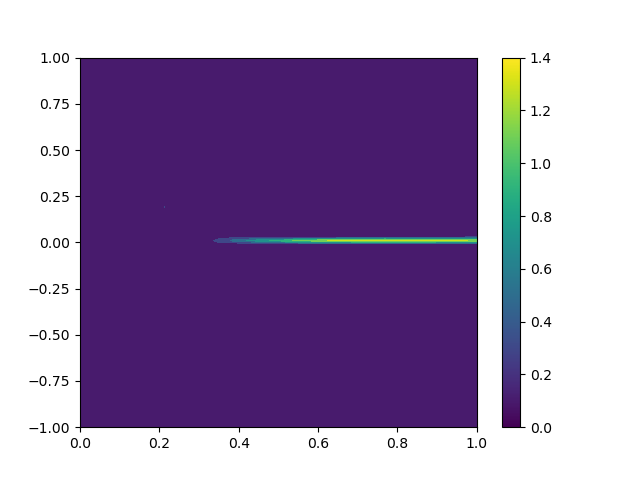}
    \caption{Absolute Error using MSE loss.}
\end{subfigure}
\begin{subfigure}{0.45\textwidth}
    \includegraphics[width=\textwidth]{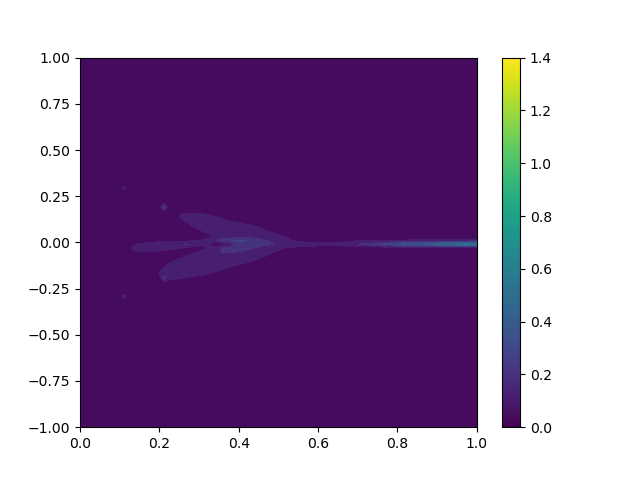}
    \caption{Absolute Error using the new loss.}
\end{subfigure}
\caption{Burger's equation solutions (time in the x-axis; space in the y-axis), (a) analytical, using both (b) classical mean loss ($\alpha=1.0$) and (c) the novel loss function ($\alpha=0.8$), along with absolute errors in both cases (d) and (e).}
\label{fig:burger_sols}
\end{figure}

The analytical and PINN solutions in time and space, along with the corresponding absolute error map are plotted for both $\alpha=0.8$ and $\alpha=1$ in figure~\ref{fig:burger_sols}.

As shown in figure~\ref{fig:burger_sols}, the $L_{\infty}$ norm of the absolute error field is reduced by nearly a factor of 2 with the proposed loss with $\alpha=0.8$ compared to the mean squared error case.

\subsection{2D linear elastic example}\label{sec:solid}

\begin{figure}[H]
    \centering
    \includegraphics[width=0.35\linewidth]{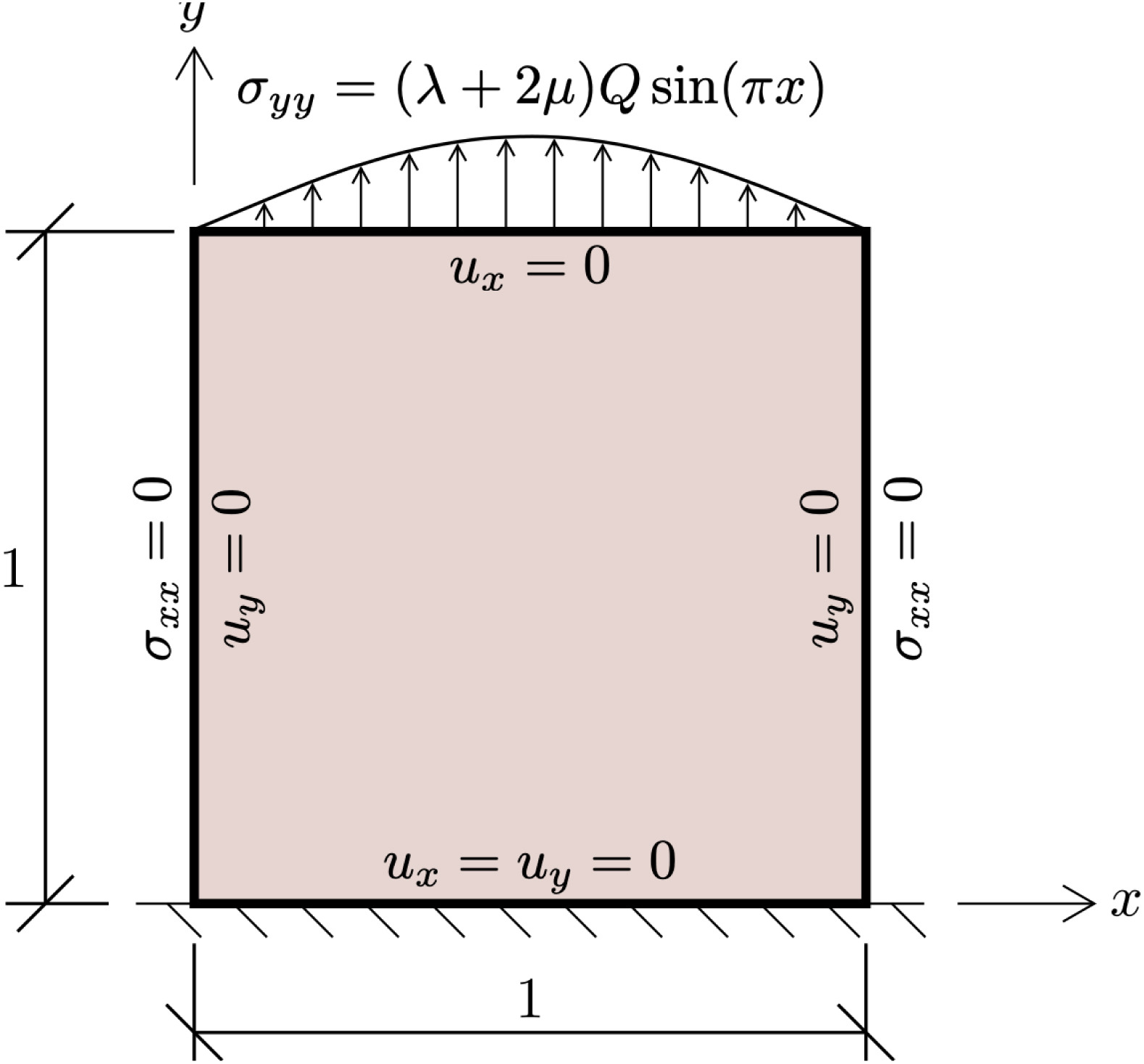}
    \caption{Solid mechanics problem setup and boundary conditions \cite{solid_pinn}.}
    \label{fig:solid_2D}
\end{figure}

We now move to a multi-dimensional example and solve a linear elastic solid mechanics problem in two dimensions as introduced by Haghighat et al. \cite{solid_pinn}. The problem definition is given in figure\ref{fig:solid_2D}.

The equations to be considered are the momentum balance, constitutive model and kinematic relationship written as:

\begin{equation}
\begin{aligned}
    &\frac{\partial \sigma_{ij}}{\partial x_j} + f_i = 0, \\[8pt]
    &\sigma_{ij} = \lambda \delta_{ij} \varepsilon_{kk} + 2\mu \varepsilon_{ij}, \\[8pt]
    &\varepsilon_{ij} = \frac{1}{2} \left( \frac{\partial u_i}{\partial x_j} + \frac{\partial u_j}{\partial x_i} \right).
\end{aligned}
\end{equation}

where in this example

\begin{equation}
\begin{aligned}
    f_x &= \lambda \left[ 4\pi^2 \cos(2\pi x) \sin(\pi y) - \pi \cos(\pi x) Q y^3 \right] 
    + \mu \left[ 9\pi^2 \cos(2\pi x) \sin(\pi y) - \pi \cos(\pi x) Q y^3 \right], \\[8pt]
    f_y &= \lambda \left[ -3 \sin(\pi x) Q y^2 + 2\pi^2 \sin(2\pi x) \cos(\pi y) \right] 
    + \mu \left[ -6 \sin(\pi x) Q y^2 + 2\pi^2 \sin(2\pi x) \cos(\pi y) + \frac{\pi^2}{4} \sin(\pi x) Q y^4 \right].
\end{aligned}
\end{equation}

where $\lambda=1$, $\mu=0.5$, and $Q=4$.

As recommended in \cite{solid_pinn}, 5 neural networks are designed to approximate the displacements and the different components of the stress tensor. Each neural network is the NN described above. In our implementation, we enforced the boundary conditions exactly using distance functions. For example, $\sigma_{xx}(x,y)$ is approximated as $x\times (1-x)\times \hat{\sigma}_{xx}(x,y)$, where $\hat{\sigma}_{xx}$ is a neural network approximation.

The solutions are obtained with 40,000 Adam iterations. Training involves 2500 randomly sampled collocation points, while the boundary conditions are automatically satisfied by construction. The computational time using MSE loss is 45 seconds vs 47 seconds using the new variance-based loss.

\begin{figure}[H]
\centering
    \begin{subfigure}[b]{0.45\textwidth}
    \centering
    \includegraphics[width=250.0pt]{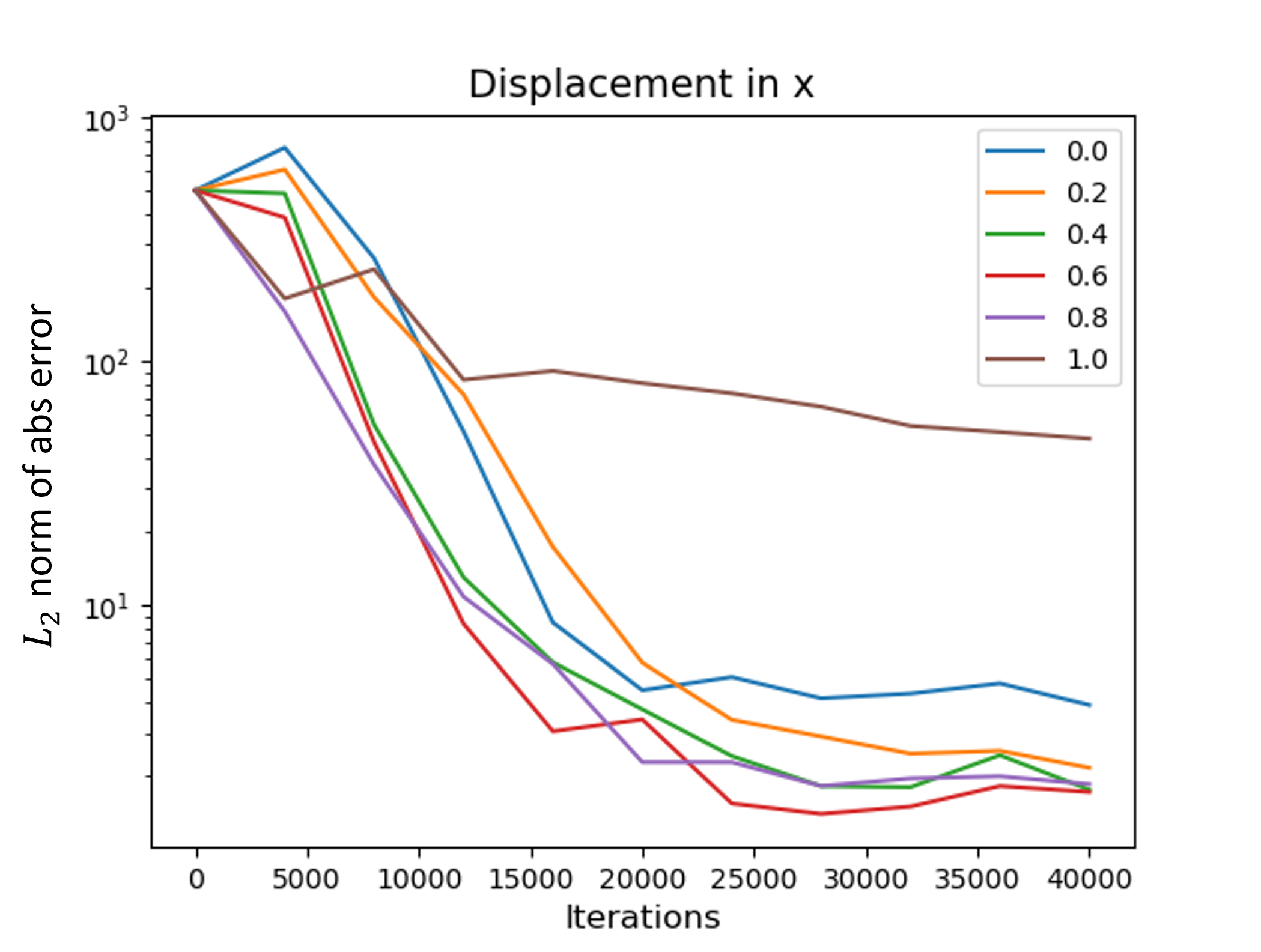}
    \end{subfigure}
    \begin{subfigure}[b]{0.45\textwidth}
    \centering
    \includegraphics[width=250.0pt]{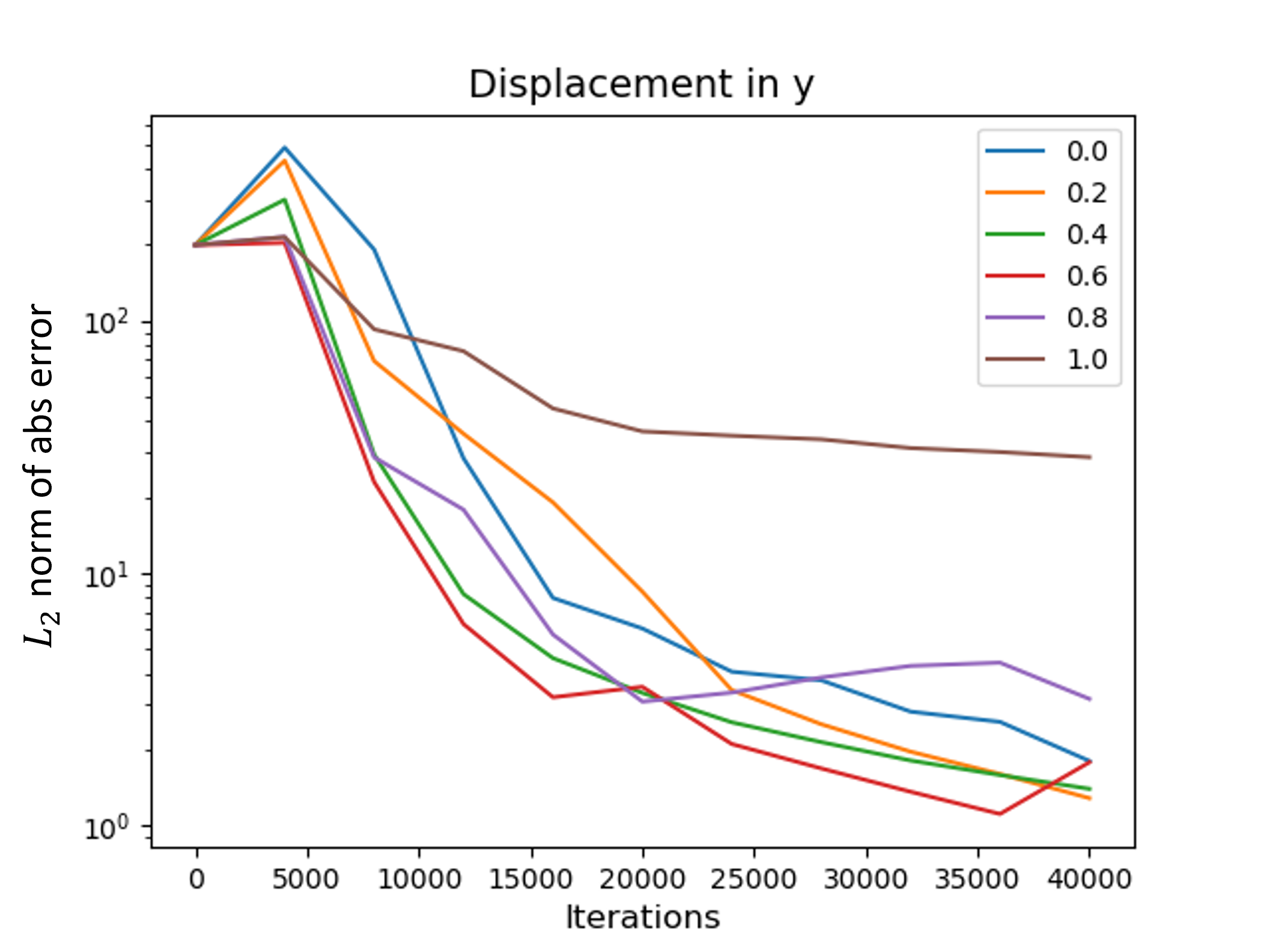}
    \end{subfigure}
\quad
\caption{The $L_2$ norm of the absolute error of x-displacement (left) and y-displacement (right), vs. the number of iterations for different values of $\alpha$ for the 2D linear elasticity example.}
\label{fig:solid_alphas}
\end{figure}

Since an analytical solution of the problem exists, we plot the $L_2$ norm of the absolute error between the predicted and analytical solutions, for both the displacements in the x and y directions, vs. the number of iterations, for different values of $\alpha$ in figure~\ref{fig:solid_alphas}.

As can be noticed from figure~\ref{fig:solid_alphas}, using the new loss function is always significantly better than using the classical mean error loss. For the displacement in the x-direction, the best solution is obtained with $\alpha=0.6$, while for the displacement in the y-direction, the best is $\alpha=0.2$ or 0.6.

\begin{figure}[H]
\centering
\begin{subfigure}{1\textwidth}
\centering
    \includegraphics[width=0.45\textwidth]{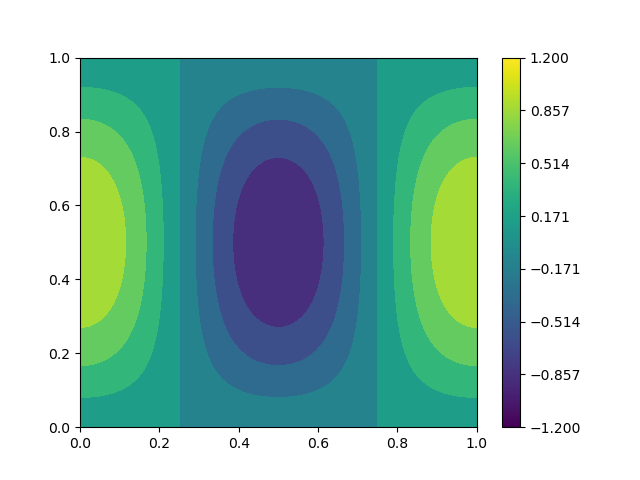}
    \caption{Analytical solution of $u_x$.}
\end{subfigure}
\begin{subfigure}{0.45\textwidth}
    \includegraphics[width=\textwidth]{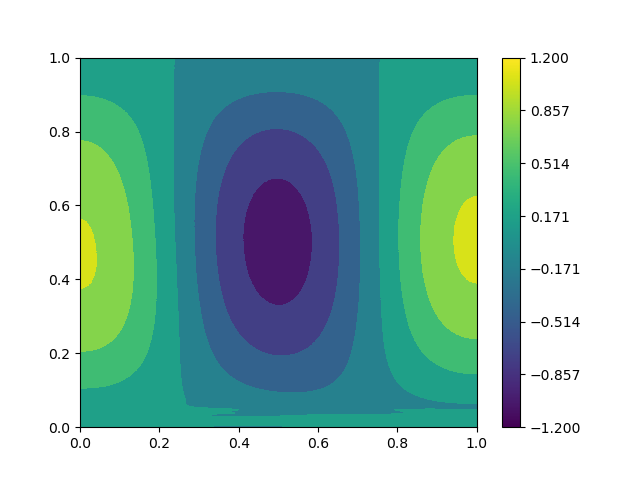}
    \caption{Solution of $u_x$ using MSE loss with $\alpha=1$.}
\end{subfigure}
\begin{subfigure}{0.45\textwidth}
    \includegraphics[width=\textwidth]{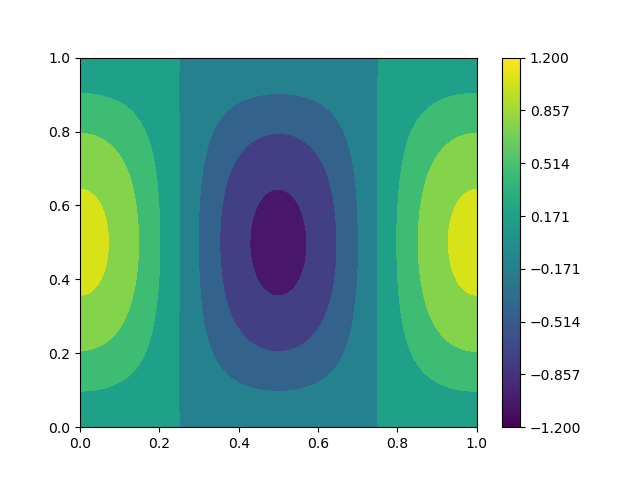}
    \caption{Solution of $u_x$ using the newly proposed loss with $\alpha=0.6$.}
\end{subfigure}
\begin{subfigure}{0.45\textwidth}
    \includegraphics[width=\textwidth]{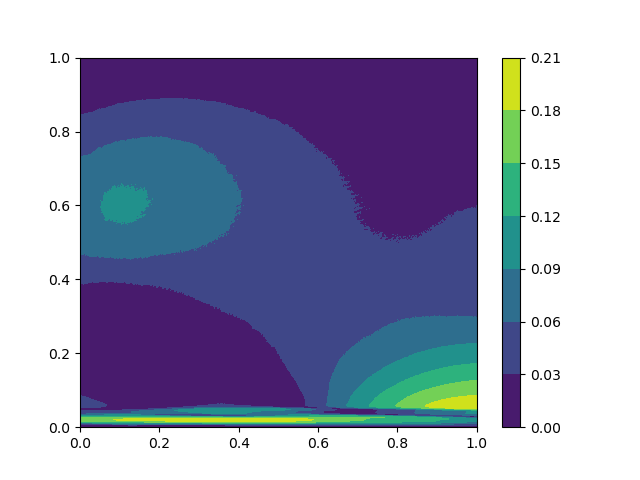}
    \caption{Absolute Error in $u_x$ using MSE loss.}
\end{subfigure}
\begin{subfigure}{0.45\textwidth}
    \includegraphics[width=\textwidth]{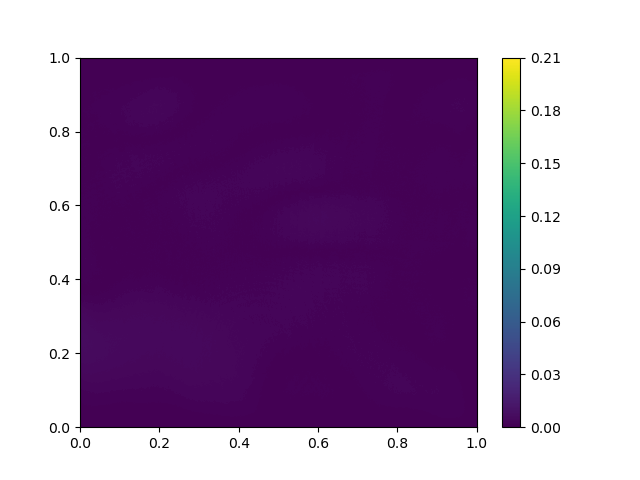}
    \caption{Absolute Error in $u_x$ using the new loss.}
\end{subfigure}
\caption{$u_x$ analytical solution along with solutions from PINN using MSE loss  ($\alpha=1$) and the new variance-based loss ($\alpha=0.6$), along with the absolute errors.}
\label{fig:solid_sol_1}
\end{figure}

\begin{figure}[H]
\centering
\begin{subfigure}{1\textwidth}
\centering
    \includegraphics[width=0.45\textwidth]{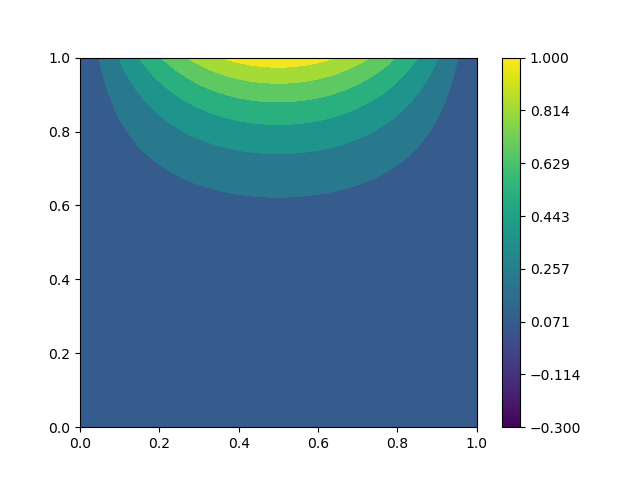}
    \caption{Analytical solution of $u_y$.}
\begin{subfigure}{0.45\textwidth}
    \includegraphics[width=\textwidth]{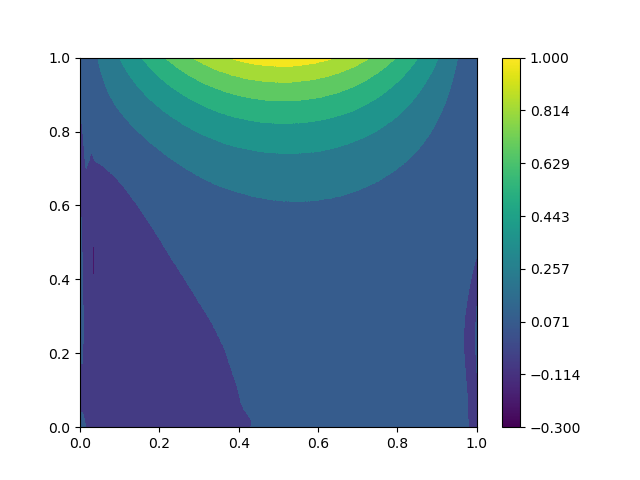}
    \caption{Solution of $u_y$ using MSE loss with $\alpha=1$.}
\end{subfigure}
\begin{subfigure}{0.45\textwidth}
    \includegraphics[width=\textwidth]{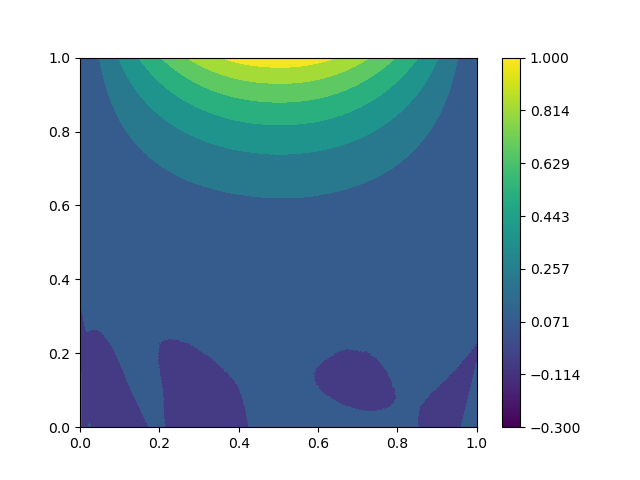}
    \caption{Solution of $u_y$ using newly proposed loss with $\alpha=0.6$.}
\end{subfigure}
\begin{subfigure}{0.45\textwidth}
    \includegraphics[width=\textwidth]{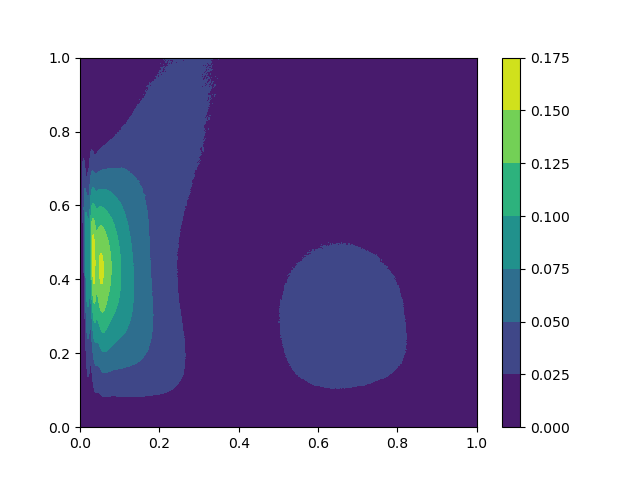}
    \caption{Absolute Error in $u_y$ using MSE loss.}
\end{subfigure}
\begin{subfigure}{0.45\textwidth}
    \includegraphics[width=\textwidth]{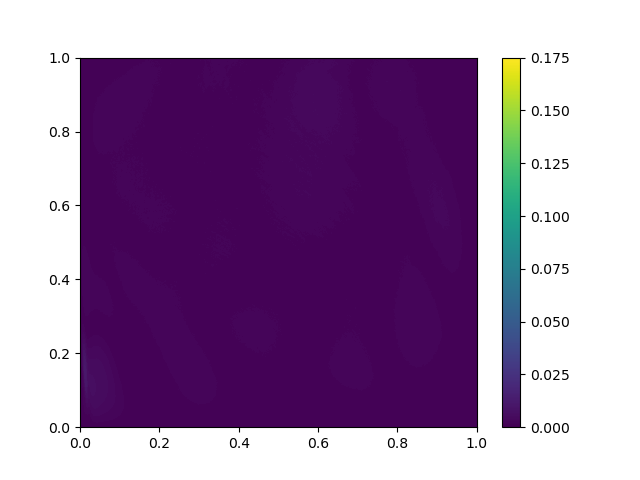}
    \caption{Absolute Error in $u_y$ using the new loss.}
\end{subfigure}
\end{subfigure}
\caption{$u_y$ analytical solution along with solutions from PINN using MSE loss  ($\alpha=1$) and the new variance-based loss ($\alpha=0.6$), along with the absolute errors.}
\label{fig:solid_sol_2}
\end{figure}

The solutions using the classical mean error loss function and the new proposed loss for $\alpha=0.6$ are shown in figures \ref{fig:solid_sol_1} and \ref{fig:solid_sol_2}.

From figures~\ref{fig:solid_sol_1} and \ref{fig:solid_sol_2}, it is clear that using the newly proposed loss function provides better solutions than using a classical mean squared error. As can be seen in for the displacement in the x direction in figure~\ref{fig:solid_sol_1}, the maximum absolute error is nearly reduced by a factor of 30. Note that the solution using MSE has a localized error region, while this is not the case for the variance-based loss.

For the displacement in the y-direction, the maximum error is reduced by a factor of 15 as shown in~figure~\ref{fig:solid_sol_2}.

\subsection{Nonlinear case in fluid mechanics}

\begin{figure}[H]
    \centering

\tikzset{every picture/.style={line width=0.75pt}} 

\begin{tikzpicture}[x=0.75pt,y=0.75pt,yscale=-1,xscale=1]

\draw  [draw opacity=0] (235.52,57.42) .. controls (235.68,57.41) and (235.84,57.41) .. (236,57.41) .. controls (252.57,57.41) and (266,70.85) .. (266,87.41) .. controls (266,87.77) and (266,88.13) .. (265.98,88.48) -- (236,87.41) -- cycle ; \draw   (235.52,57.42) .. controls (235.68,57.41) and (235.84,57.41) .. (236,57.41) .. controls (252.57,57.41) and (266,70.85) .. (266,87.41) .. controls (266,87.77) and (266,88.13) .. (265.98,88.48) ;  
\draw    (308.38,87.34) -- (307.5,159.36) ;
\draw    (265.98,88.48) -- (265.5,159.36) ;
\draw  [draw opacity=0] (235.5,23.46) .. controls (235.67,23.45) and (235.84,23.45) .. (236,23.45) .. controls (275.94,23.45) and (308.33,52.05) .. (308.38,87.34) -- (236,87.41) -- cycle ; \draw   (235.5,23.46) .. controls (235.67,23.45) and (235.84,23.45) .. (236,23.45) .. controls (275.94,23.45) and (308.33,52.05) .. (308.38,87.34) ;  
\draw    (235.5,23.46) -- (235.52,57.42) ;
\draw    (265.5,159.36) -- (307.5,159.36) ;
\draw  [draw opacity=0] (267.51,209.25) .. controls (267.45,208.01) and (267.43,206.75) .. (267.46,205.48) .. controls (267.85,187.02) and (277.38,172.25) .. (288.76,172.49) .. controls (300.13,172.73) and (309.03,187.89) .. (308.64,206.35) .. controls (308.63,207.05) and (308.6,207.76) .. (308.56,208.45) -- (288.05,205.91) -- cycle ; \draw   (267.51,209.25) .. controls (267.45,208.01) and (267.43,206.75) .. (267.46,205.48) .. controls (267.85,187.02) and (277.38,172.25) .. (288.76,172.49) .. controls (300.13,172.73) and (309.03,187.89) .. (308.64,206.35) .. controls (308.63,207.05) and (308.6,207.76) .. (308.56,208.45) ;  
\draw    (273.67,208.22) -- (273.67,186.22) ;
\draw [shift={(273.67,184.22)}, rotate = 90] [color={rgb, 255:red, 0; green, 0; blue, 0 }  ][line width=0.75]    (10.93,-3.29) .. controls (6.95,-1.4) and (3.31,-0.3) .. (0,0) .. controls (3.31,0.3) and (6.95,1.4) .. (10.93,3.29)   ;
\draw    (288.05,205.91) -- (287.69,173.22) ;
\draw [shift={(287.67,171.22)}, rotate = 89.36] [color={rgb, 255:red, 0; green, 0; blue, 0 }  ][line width=0.75]    (10.93,-3.29) .. controls (6.95,-1.4) and (3.31,-0.3) .. (0,0) .. controls (3.31,0.3) and (6.95,1.4) .. (10.93,3.29)   ;
\draw    (301.67,207.22) -- (301.67,184.22) ;
\draw [shift={(301.67,182.22)}, rotate = 90] [color={rgb, 255:red, 0; green, 0; blue, 0 }  ][line width=0.75]    (10.93,-3.29) .. controls (6.95,-1.4) and (3.31,-0.3) .. (0,0) .. controls (3.31,0.3) and (6.95,1.4) .. (10.93,3.29)   ;

\draw (228,209) node [anchor=north west][inner sep=0.75pt]   [align=left] {Inlet velocity profile};
\draw (133,29) node [anchor=north west][inner sep=0.75pt]    {$P=0\ ( outlet)$};
\draw (316.95,40.45) node [anchor=north west][inner sep=0.75pt]  [rotate=-68.23] [align=left] {No-slip wall};

\end{tikzpicture}
    \caption{Fluid mechanics problem geometry and boundary conditions.}
    \label{fig:fluid_problem}
\end{figure}
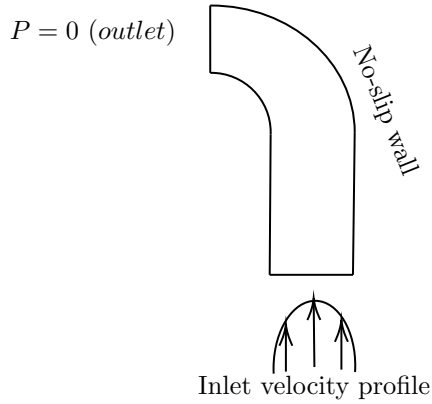

In this example, we consider a nonlinear problem, namely the steady state Navier-Stokes equations in 2D. The description of the problem is provided in figure~\ref{fig:fluid_problem}. An inlet velocity is prescribed with a parabolic profile in the direction normal to the boundary,  and zero pressure is prescribed on the outlet boundary. The rest of the boundary (the wall) has a no slip boundary condition.

The PDE system to be solved is:

\begin{equation}
    \rho \left( u \frac{\partial u}{\partial x} + v \frac{\partial u}{\partial y} \right) = -\frac{\partial p}{\partial x} + \mu \left( \frac{\partial^2 u}{\partial x^2} + \frac{\partial^2 u}{\partial y^2} \right),
\end{equation}
\begin{equation}
    \rho \left( u \frac{\partial v}{\partial x} + v \frac{\partial v}{\partial y} \right) = -\frac{\partial p}{\partial y} + \mu \left( \frac{\partial^2 v}{\partial x^2} + \frac{\partial^2 v}{\partial y^2} \right), \\
\end{equation}

\begin{equation}
    \frac{\partial u}{\partial x} + \frac{\partial v}{\partial y} = 0.
\end{equation}

Here, \(u, v\) are the velocity components in the \(x\)- and \(y\)-directions, respectively; \(\rho\) is the fluid density, \(p\) is the pressure, and \(\mu\) is the dynamic viscosity.

\begin{figure}[H]
\centering
\begin{subfigure}{1\textwidth}
\centering
    \includegraphics[width=0.45\textwidth]{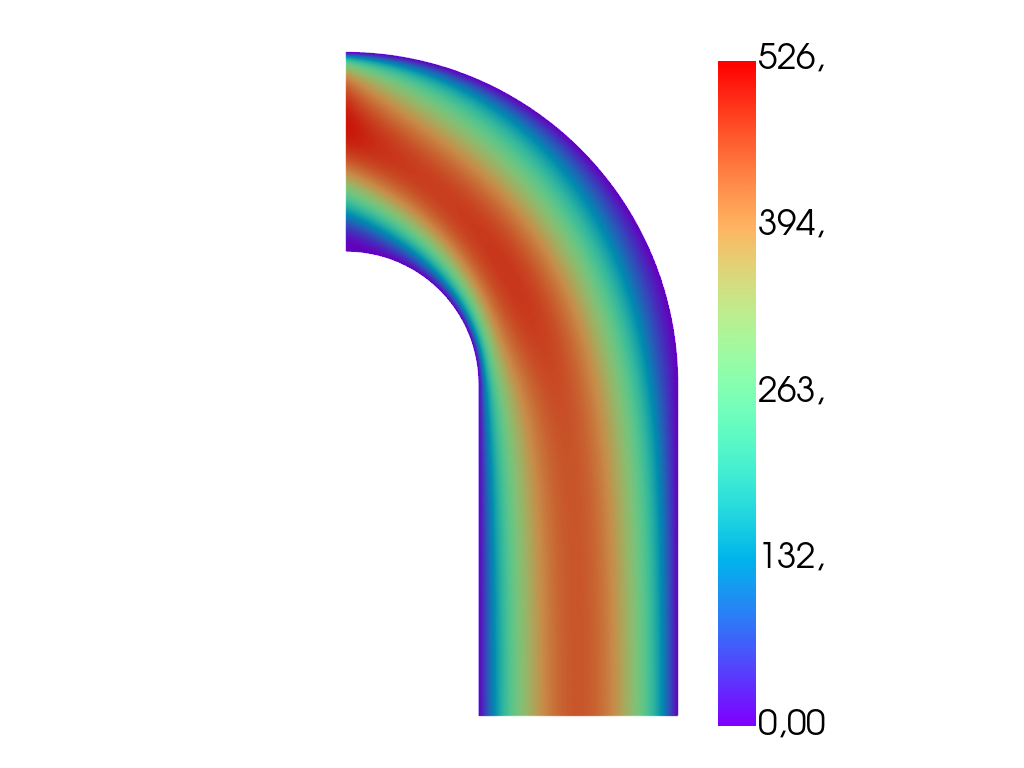}
    \caption{Velocity magnitude from CFD (reference solution).}
    \label{fig:vel_cfd}
\end{subfigure}
\begin{subfigure}{0.45\textwidth}
    \includegraphics[width=\textwidth]{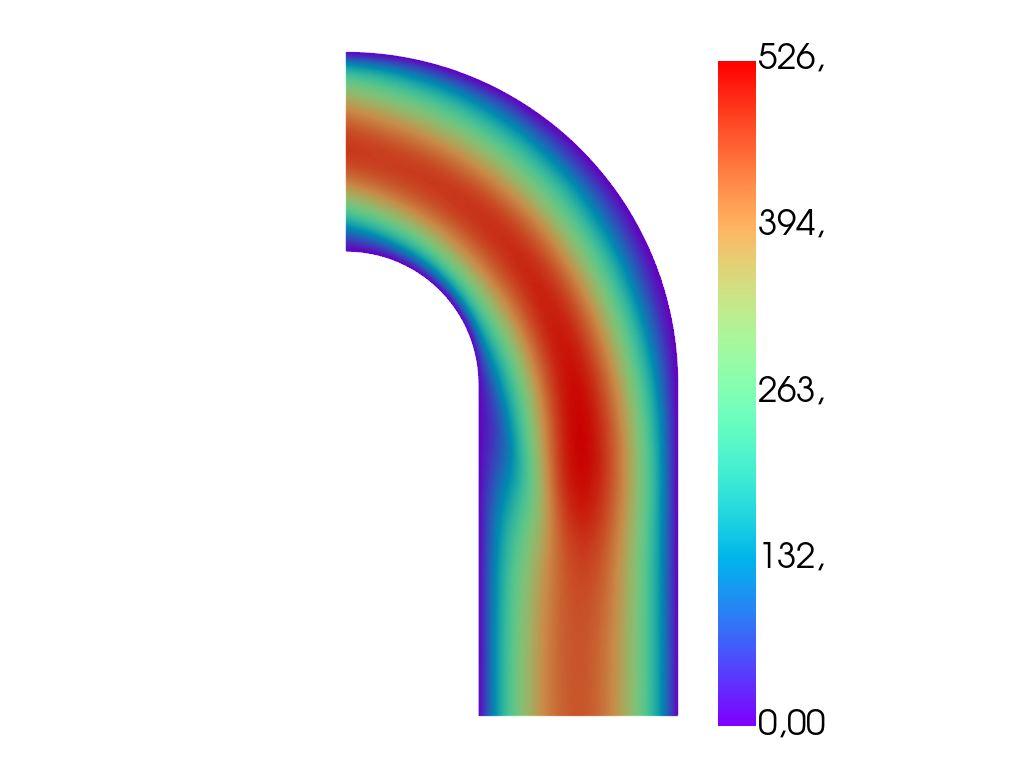}
    \caption{Velocity magnitude obtained using MSE in PINN.}
    \label{fig:vel_mse}
\end{subfigure}
\begin{subfigure}{0.45\textwidth}
    \includegraphics[width=\textwidth]{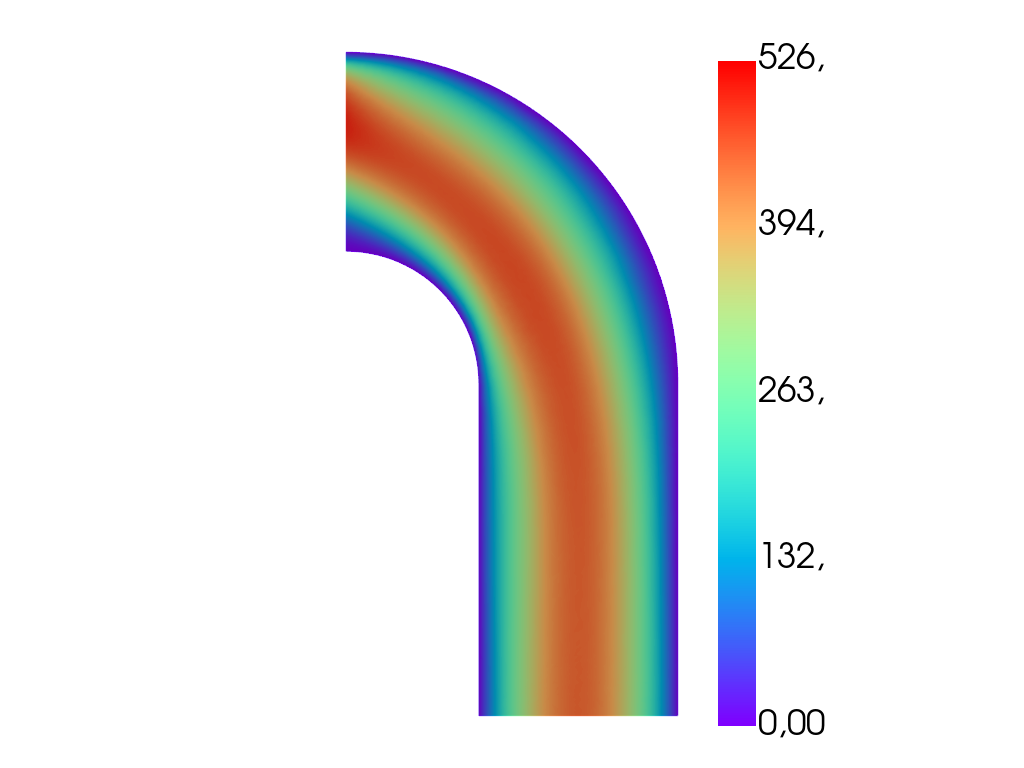}
    \caption{Velocity magnitude obtained using the new loss function in PINN.}
    \label{fig:vel_std}
\end{subfigure}
\begin{subfigure}{0.45\textwidth}
    \includegraphics[width=\textwidth]{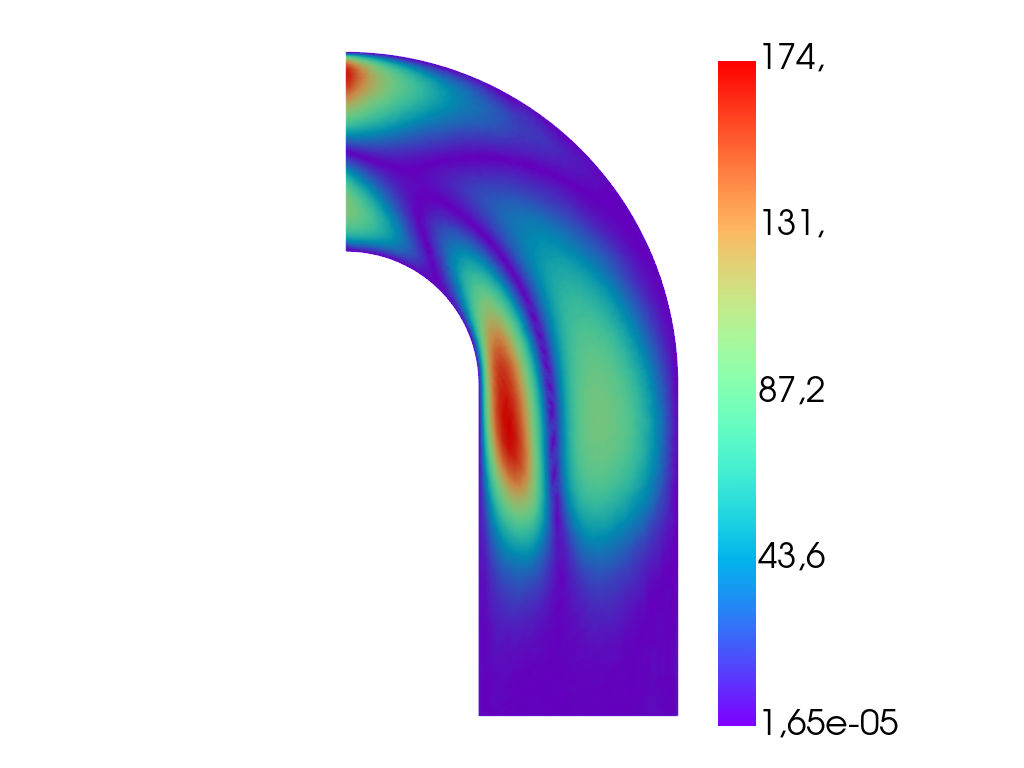}
    \caption{Absolute error in velocity magnitude using MSE loss.}
    \label{fig:vel_mse}
\end{subfigure}
\begin{subfigure}{0.45\textwidth}
    \includegraphics[width=\textwidth]{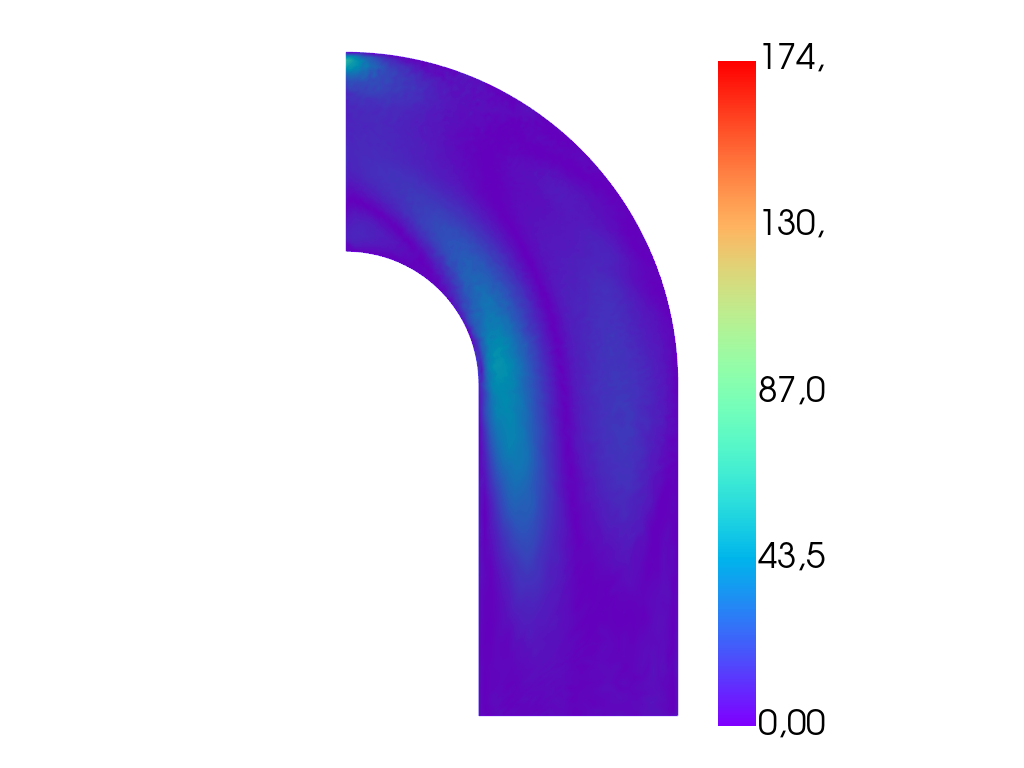}
    \caption{Absolute error in velocity magnitude using the new loss.}
    \label{fig:vel_std}
\end{subfigure}
\caption{Velocity magnitude using Fenics as a reference solution along with solutions from PINN using MSE loss  ($\alpha=1$) and the new variance-based loss ($\alpha=0.8$), along with the absolute errors.}
\label{fig:velocities_mag}
\end{figure}

In this example the maximum velocity is chosen to be 450 m/s, with a density and viscosity of unity. These choices given the geometry at hand lead to a non-parabolic velocity profile at the outlet after the curvature as observed from the reference solution in figure~\ref{fig:vel_cfd}.

The no-slip boundary condition is enforced automatically with an approximation of the distance function using a separate neural network which was previously trained before PINN training. Moreover, in order to avoid differentiating the neural networks twice, the derivatives of the velocity components are approximated with separate neural networks and extra terms are added in the loss function to enforce their relationship with the velocity networks.

\begin{figure}[H]
\centering
    \begin{subfigure}[b]{0.45\textwidth}
    \centering
    \includegraphics[width=250.0pt]{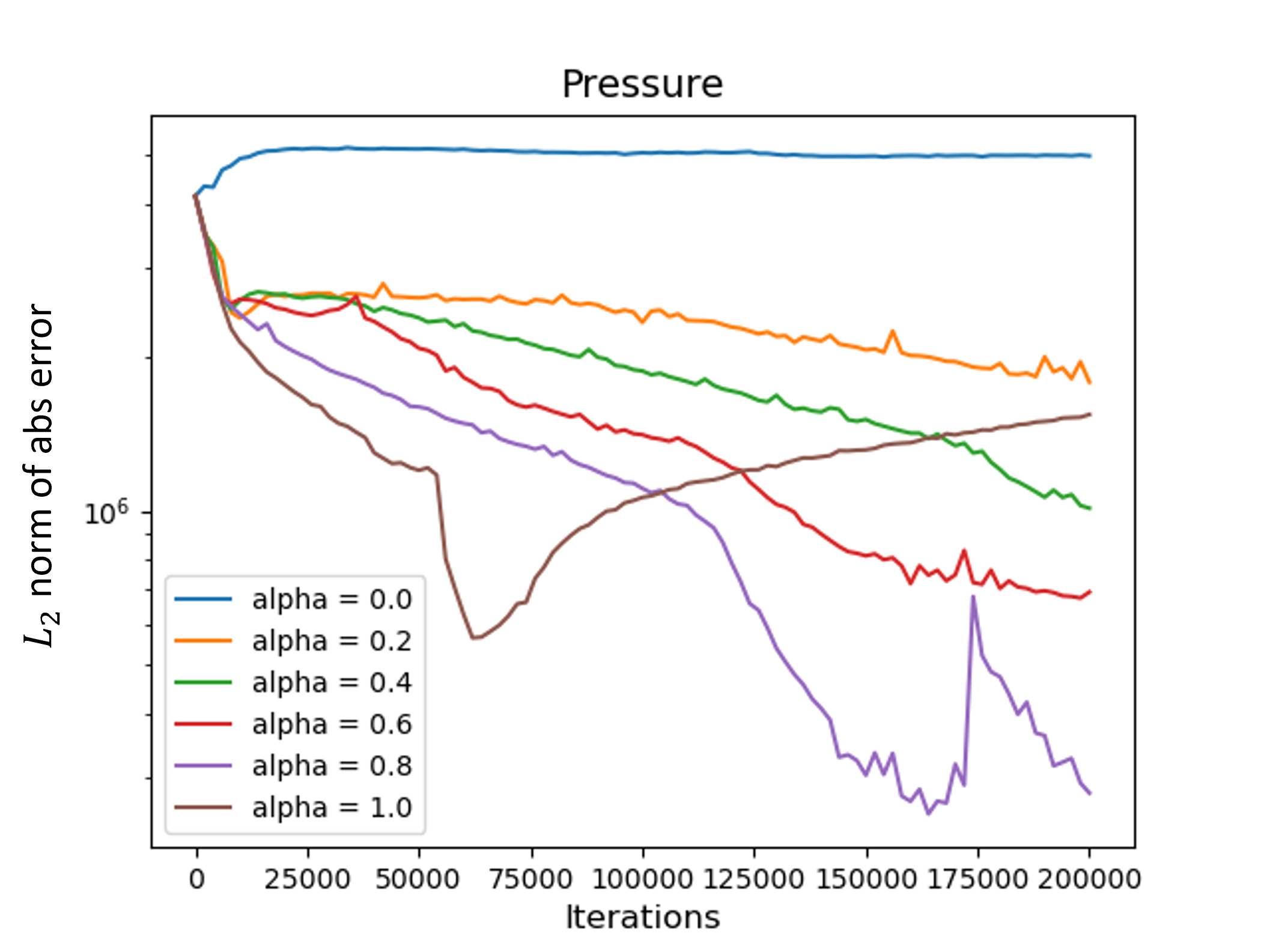}
    \end{subfigure}
    \begin{subfigure}[b]{0.45\textwidth}
    \centering
    \includegraphics[width=250.0pt]{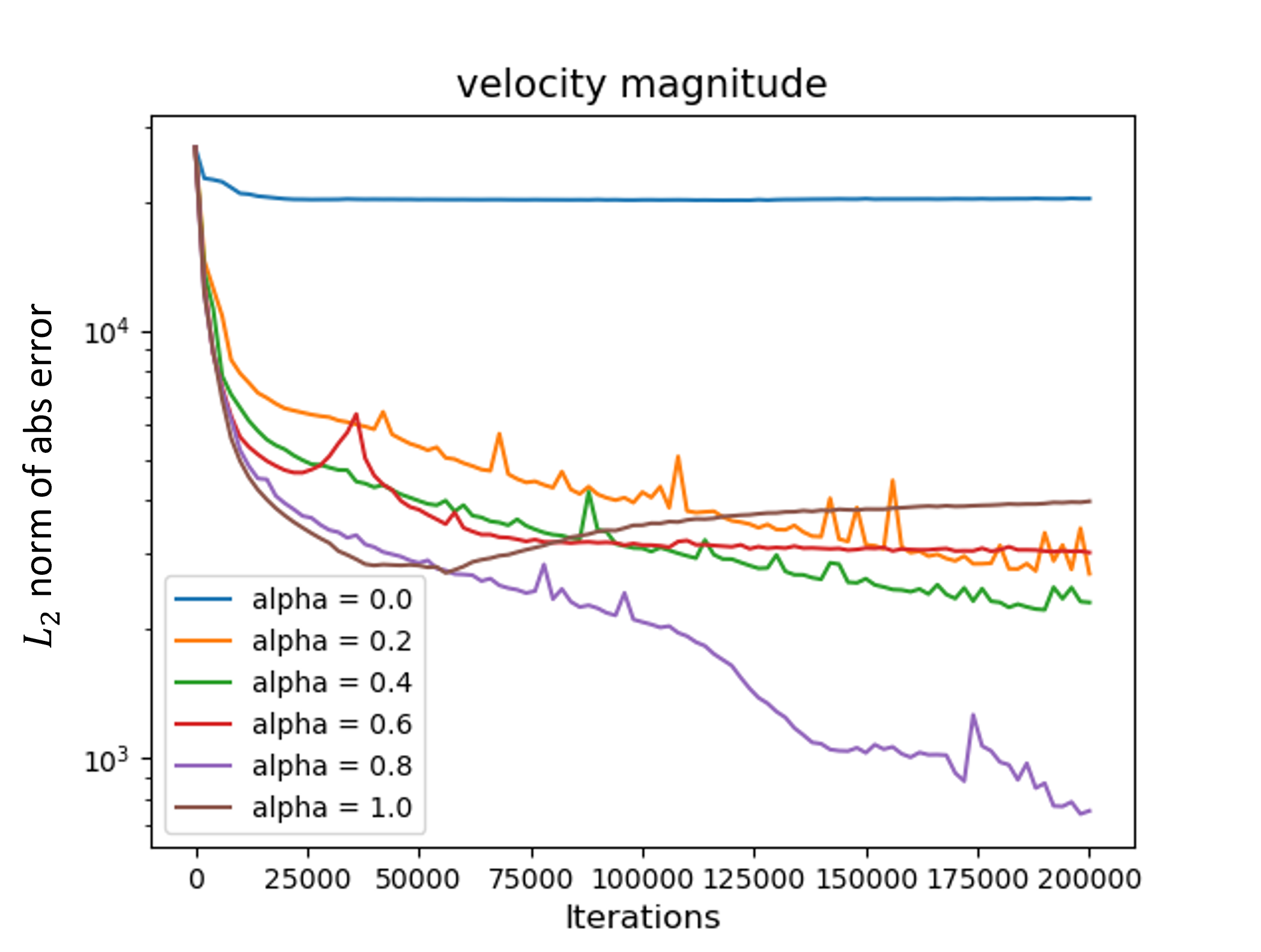}
    \end{subfigure}
\quad
\caption{On the left, the $L_2$ norm of the absolute error of the pressure vs. the number of iterations for different values of $\alpha$ for the 2D Navier-Stokes example. On the right, the $L_2$ norm of the absolute error of velocity magnitude vs. the number of iterations for different values of $\alpha$.}
\label{fig:cfd_alphas}
\end{figure}

In this example, 200,000 Adam iterations are employed. The computational time is nearly 15 minutes for MSE loss vs 16 minutes for our proposed loss. Training involves 10,176 collocation points, generated from meshing the domain with triangular elements, which also serve as input for the CFD simulation (P1–P1 finite elements with Fenics). The $L_2$ norm of the absolute error between the PINN and CFD solutions is plotted in figure~\ref{fig:cfd_alphas} for different $\alpha$ values.

It can be noticed from figure~\ref{fig:cfd_alphas} that the mean squared loss ($\alpha=1$) tends to diverge from the reference solution as opposed to the converged solutions with the new proposed loss, with the best solution by choosing $\alpha=0.8$.

\begin{figure}[H]
\centering
\begin{subfigure}{1\textwidth}
\centering
    \includegraphics[width=0.45\textwidth]{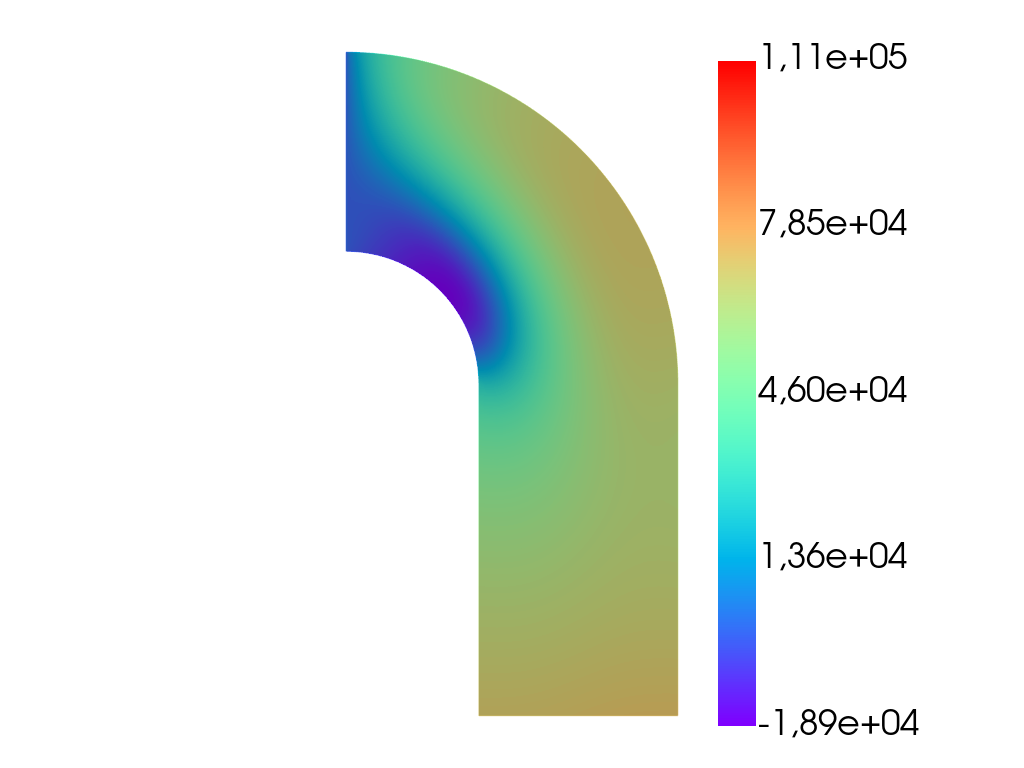}
    \caption{Pressure field from CFD (reference solution).}
    \label{fig:first}
\end{subfigure}
\begin{subfigure}{0.45\textwidth}
    \includegraphics[width=\textwidth]{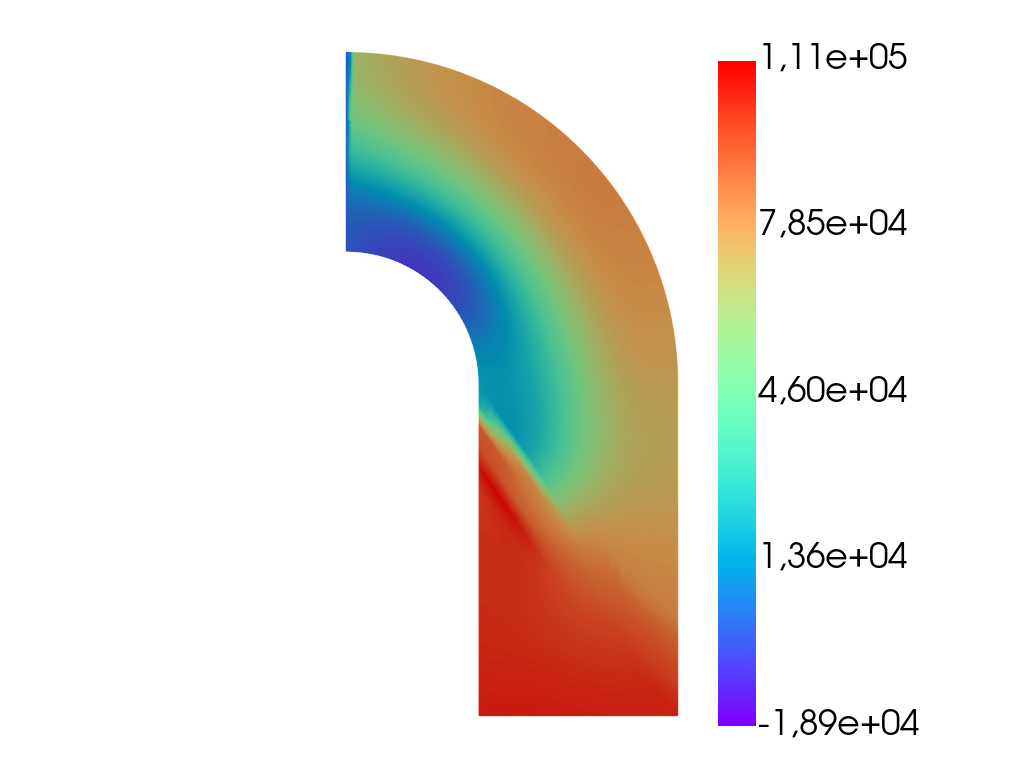}
    \caption{Pressure field obtained using PINN with MSE loss function.}
    \label{fig:second}
\end{subfigure}
\begin{subfigure}{0.45\textwidth}
    \includegraphics[width=\textwidth]{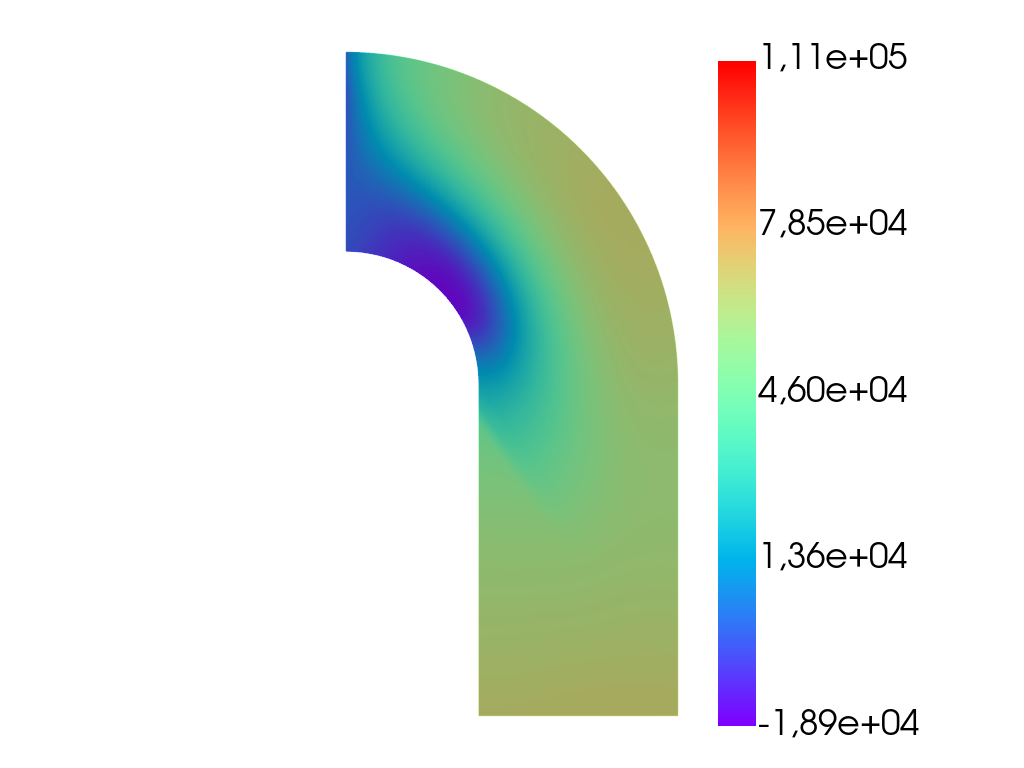}
    \caption{Pressure field obtained using PINN with the new proposed loss.}
    \label{fig:third}
\end{subfigure}
\begin{subfigure}{0.45\textwidth}
    \includegraphics[width=\textwidth]{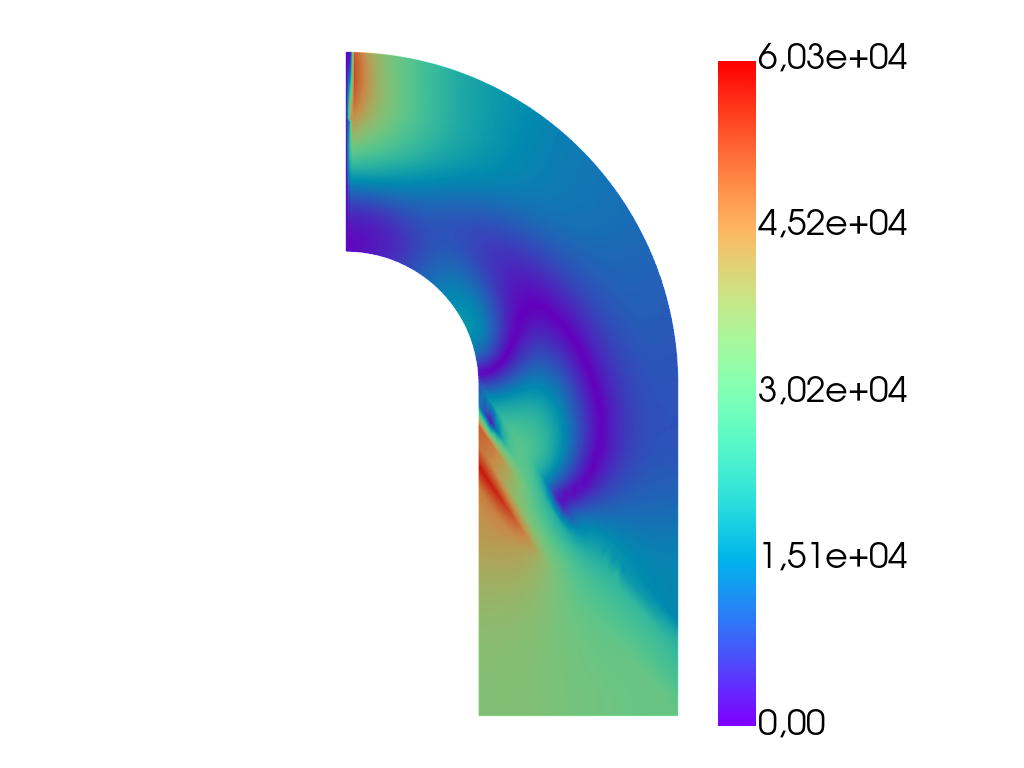}
    \caption{Absolute error in pressure field using PINN with MSE loss.}
    \label{fig:second}
\end{subfigure}
\begin{subfigure}{0.45\textwidth}
    \includegraphics[width=\textwidth]{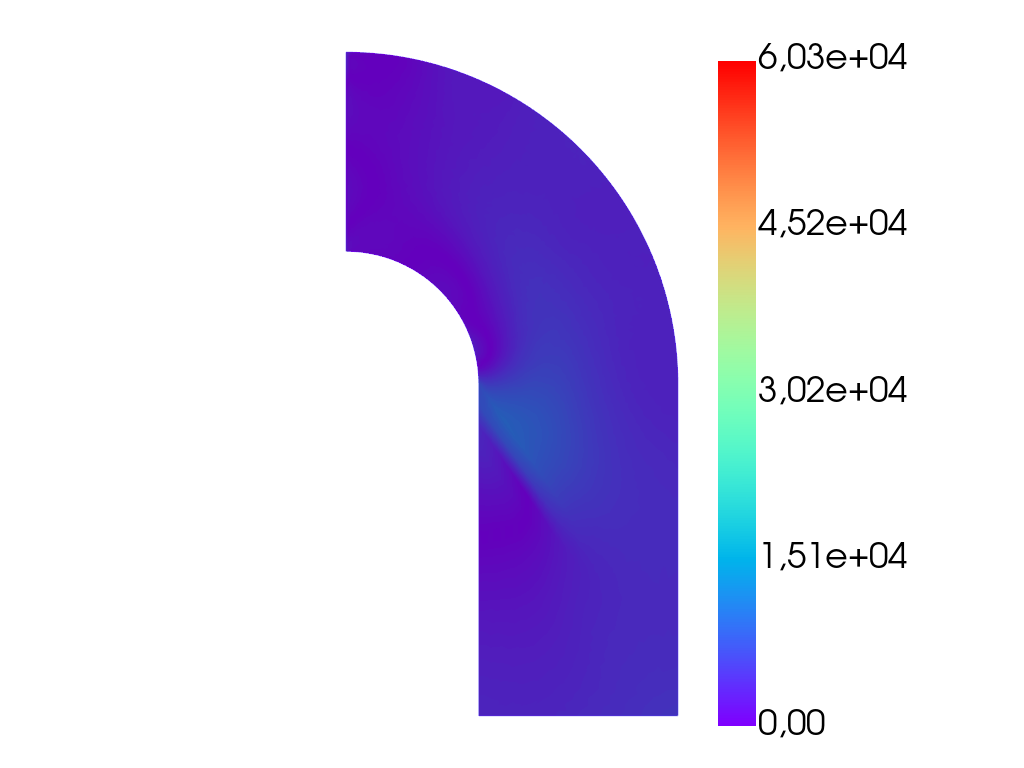}
    \caption{Absolute error in pressure field using PINN with the new loss.}
    \label{fig:third}
\end{subfigure}
\caption{Pressure fields using CFD (reference solution), PINN with MSE loss  ($\alpha=1$) and PINN with the new variance-based loss  ($\alpha=0.8$), along with the absolute errors.}
\label{fig:pressures}
\end{figure}

The results are shown in figures~\ref{fig:velocities_mag} and \ref{fig:pressures} for the mean squared loss and using variance-based loss with $\alpha=0.8$.

As can be seen from figure~\ref{fig:velocities_mag}, the velocity magnitude reference solution has an outlet velocity that has a shifted maximum towards the top due to centrifugal forces. This behavior is captured when using the variance-based loss, but not when using the MSE loss. It can be clearly noticed that the maximum error value is reduced by a factor of more than 3 using the new proposed loss.

Regarding pressure (figure~\ref{fig:pressures}), it is noticed the the MSE loss solution is not smooth and does not match the reference well. By contrast, the new loss function solution is closer to the reference solution and smooth. Moreover, the maximum error in the pressure estimation with the variance-based loss is one order of magnitude less than with the MSE loss.

\section{Discussion and conclusion}\label{sec:discuss}

\subsection{Computational cost and implementation}

The computational cost incurred by adding this extra term in the loss function is almost negligible, since the operations needed to calculate the standard deviation or its derivative in the computational graph are basic.

Implementation-wise, the new term is quite straightforward to add to the loss function. It typically consists of changing one line in the code. The term can be added in TensorFlow using the \texttt{tf.math.reduce\_std} function. In PyTorch, it can be introduced using the \texttt{torch.std} function.

\subsection{Comparison with Gradient-enhanced PINN and Huber loss}

\begin{figure}[H]
\centering
    \begin{subfigure}[b]{0.45\textwidth}
    \centering
    \includegraphics[width=250.0pt]{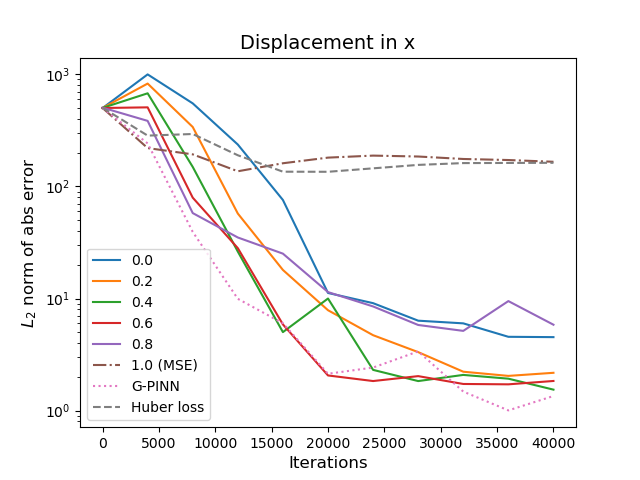}
    \end{subfigure}
    \begin{subfigure}[b]{0.45\textwidth}
    \centering
    \includegraphics[width=250.0pt]{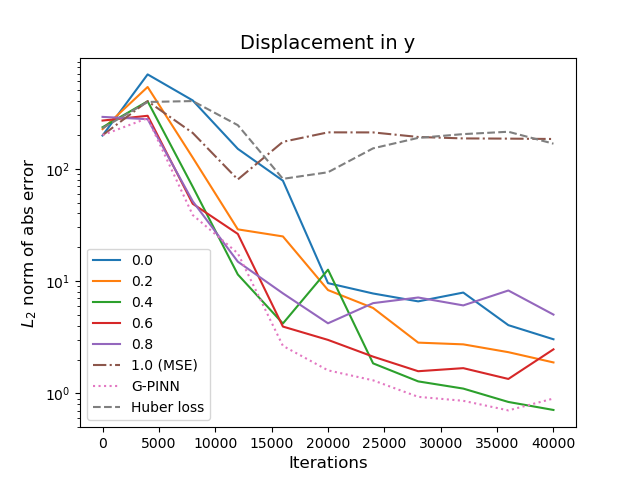}
    \end{subfigure}
\quad
\caption{On the left, the $L_2$ norm of the absolute error of x-displacement vs. the number of iterations for different values of $\alpha$ for the 2D linear elasticity example along with gPINNs and Huber loss solutions. On the right, the $L_2$ norm of the absolute error of x-displacement vs. the number of iterations for different values of $\alpha$ along with gPINNs and Huber loss solutions.}
\label{fig:gpinn}
\end{figure}

A numerical comparison is performed using Gradient-enhanced PINN (gPINNs) \cite{gpinn} and Huber loss along with our proposed loss for the solid mechanics example introduced in section~\ref{sec:solid}. The error plots are shown in figure~\ref{fig:gpinn}.

For the Huber loss case, we chose $\delta=1$, where $\delta$ controls the steepness of the Huber curve. From the figure, we can see that the Huber loss has almost no improvement compared with the MSE loss. ($\alpha=1$). 

gPINNs has a similar effect as using our new proposed loss function. However, it provides a marginal improvement to the new loss at the expense of a high computational cost. The cost of gPINN for this problem is greater than 3$\times$ than that of the original PINN, whereas our proposed modification has negligible computational overhead. Moreover, gPINNs is more difficult to implement, while our proposed regularization is straightforward.

\subsection{Effect on collocation points' density}

\begin{figure}[H]
\centering
    \begin{subfigure}[b]{0.45\textwidth}
    \centering
    \includegraphics[width=250.0pt]{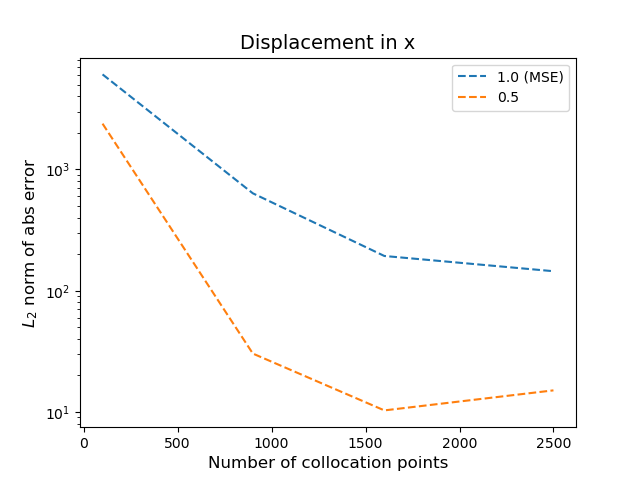}
    \end{subfigure}
    \begin{subfigure}[b]{0.45\textwidth}
    \centering
    \includegraphics[width=250.0pt]{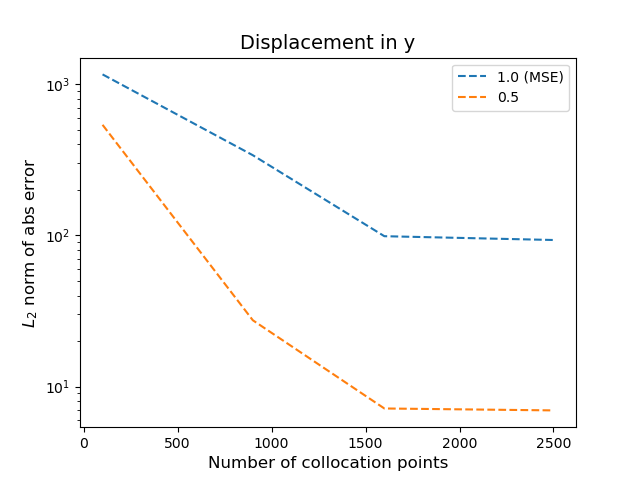}
    \end{subfigure}
\quad
\caption{On the left, the $L_2$ norm of the absolute error of x-displacement vs. the number of collocation points for $\alpha$ values $1.0$ and $0.5$, for the 2D linear elasticity example. On the left, the $L_2$ norm of the absolute error of y-displacement vs. the number of collocation points for $\alpha$ values $1.0$ and $0.5$, for the 2D linear elasticity example.}
\label{fig:cols}
\end{figure}

In this section, we study the effect of the added term in the loss function on the collocation point's density needed to reach a comparable accuracy using MSE loss. In figure~\ref{fig:cols}, the $L_2$ norm of the error is plotted for two cases of $\alpha$ (1.0, 0.5) vs. the number of collocation points.

From figure~\ref{fig:cols}, the influence of the additional term in the loss function is evident. In particular, for $\alpha = 0.5$, the model achieves a given level of accuracy with substantially fewer collocation points compared to the standard MSE loss ($\alpha = 1.0$). This suggests that the modified loss function improves the efficiency of the training process by enforcing the governing equations more effectively, thus reducing the reliance on dense spatial sampling.

Moreover, the improved performance indicates that the added term enhances the network's ability to capture the underlying physical behavior with fewer constraints. This can be especially beneficial in high-dimensional or computationally expensive problems, where training on a large number of collocation points is impractical or computationally expensive.

\subsection{Conclusion}
In this article, we have proposed a new loss function combining MSE and variance to regularize PINNs. In a variety of examples, the PINN solution clearly improves when applying the new variance-based loss function. This modification effectively reduces the maximum error and minimizes localized error regions, leading to an overall enhancement in solution behavior. The new loss function is also more efficient than the MSE loss, requiring less collocation points: this could be a critical point in higher dimension problems. Moreover, it appears better in accuracy or efficiency than two other existing enhanced methods, the Huber loss and gPINNs. 

In most of our cases, we empirically found that $\alpha$ = 0.8 often yields the best results: a theoretical justification remains an open question.

Only the Adam optimizer was employed in all the cases: in the Burgers and linear elasticity examples, this optimizer is enough to obtain good accuracy. For the Navier-Stokes example, Adam did not provide the best solution. In future work, we will explore whether a more powerful optimizer can be applied along with Adam such as BFGS to improve the solution. But for the purpose of testing our new loss function, Adam optimizer was enough to demonstrate the usefulness of our proposal.

In conclusion, the added term to the loss function is simple to implement, and the low computational overhead is largely outweighed by the solution improvement. 

\section*{Funding Acknowledgment}

We acknowledge funding from the European Research Council (ERC) under the European Union’s Horizon 2020 Research and Innovation Program (Grant Agreement No. 864313).

\bibliography{biblio}

\end{document}